\documentclass[12pt]{amsart}
\usepackage{fourier,microtype}
\usepackage{geometry}
\usepackage{amsfonts, amsmath, amsthm}
\usepackage{graphicx}
\usepackage{epsfig}
\usepackage{wrapfig}

\newtheorem{pro}{Proposition}[section]
\newtheorem{thm}[pro]{Theorem}
\newtheorem{lem}[pro]{Lemma}

\newtheorem{props}[pro]{Properties}

\newtheorem{conds}[pro]{Conditions}
\newtheorem{clm}{Claim}

\newtheorem{question}[pro]{Question}
\newtheorem{example}[pro]{Example}

\newtheorem{rmkk}[pro]{Remark}

\newtheorem{inv}{Invariant}

\theoremstyle{definition}
\newtheorem{dfn}[pro]{Definition}
\newtheorem{dfns}[pro]{Definitions}

\theoremstyle{remark}
\newtheorem*{rmk}{Remark}
\newtheorem*{rmks}{Remarks}

\newcommand{\bfP}{{\bf P}}
\newcommand{\bfp}{{\bf p}}
\newcommand{\bfY}{{\bf Y}}
\newcommand{\bfy}{{\bf y}}

\newcommand{\tpe}{t_0+\epsilon}
\newcommand{\scc}{simple closed curve}
\newcommand{\s}{\Sigma}

\newcommand{\gst}{f(\Sigma_t)}
\newcommand{\sfs}{Seifert fibered space}

\newcommand{\ie}{{\it i.e.}}

\newcommand{\ab}{{\it ``above"}}
\newcommand{\bel}{{\it ``below"}}
\newcommand{\eg}{{\it e.g.}}
\newcommand{\bbb}{\mathbb}
\newcommand{\del}{\partial}
\newcommand{\hh}{Heegaard}
\newcommand{\assu}{an irreducible, orientable, a-toroidal, non-Seifert fibered}
\newcommand{\hgt}{of Heegaard genus at least two}
\newcommand{\hhs}{Heegaard surface}
\newcommand{\hhf}{Heegaard function}

\newcommand{\ka}{\textsf{K}}
\newcommand{\ca}{\textsf{C}}
\newcommand{\sca}{\mbox{\rm sing}(\ca)}

\newcommand{\mst}{maximal solid torus}
\newcommand{\msts}{maximal solid tori}
\newcommand{\st}{solid torus}

\newcommand{\op}{orientation preserving}

\newcommand{\bddup}{8 g(\s) - 7}
\newcommand{\bddupone}{8 g(\s) - 6}
\newcommand{\bdddown}{4g(\s) - 3}

\title[Invariant Heegaard Surfaces in Manifolds with Involutions]{Invariant
  {H}eegaard Surfaces \\ in Manifolds with Involutions \\ 
  and the {H}eegaard genus of double covers}

\date{December 25, 2009}
\address{Department of mathematical Sciences, University of
Arkansas, Fayetteville, AR 72701}
\email{yoav@math.uark.edu}
\author{Yo'av Rieck}
\address{Department of Mathematics and Statistics, University of
Melbourne, Melbourne, Parkville 3010, Australia}
\email{rubin@ms.unimelb.edu.au}
\author{J Hyam Rubinstein}
\thanks{First named author supported in part by JSPS grant P00024 and  
 The 21st Century COE Program ``Constitution of wide-angle mathematical basis 
 focused on knots'' (Project Leader¡§Akio Kawauchi)} 

\begin{document}

\begin{abstract}

Let $M$ be a 3-manifold admitting a strongly irreducible Heegaard surface $\s$
and $f:M\to M$ an involution.  We construct an invariant Heegaard surface for
$M$ of genus at most $\bddup$.  As a consequence, given a
(possibly branched)  double cover $\pi:M \to N$ we obtain the following bound
on the Heegaard genus of $N$:

$$g(N) \leq \bdddown. $$

\noindent We also get a bound on the complexity of the branch set in
terms of $g(\s)$.  If we assume that $M$ is non-Haken, by Casson and Gordon
\cite{casson-gordon} we may replace $g(\s)$ by $g(M)$ in all the statements
above.
\end{abstract}

\maketitle

\section{Statements of results}

We study the behavior of \hhs s under (possibly branched) double covers $\pi:M
\to N$.   It is easy to lift any \hhs\ of $N$ to a \hhs\ of $M$ and see that
the \hh\ genus of the cover is bounded above: $g(M) \leq 2g(N) + b -1$.
Here $g(\cdot)$ denotes the genus of a surface or the Heegaard genus of a
3-manifold and $b$ is the bridge index of the branch set with respect to a
minimal genus Heegaard surface for $N$.  (This upper bound easily generalizes
to any $p$-fold cover $\pi:M \to N$, provided that the branch set is a
1-manifold: $g(M) \leq p g(N) + (p-1) b -1$; see, for example,
\cite[Chapter~11]{burde-z}.)

For the converse we need a strongly irreducible \hhs\ for $M$, say $\s$.
Since any double cover is regular, it is given as the quotient under an
involution $f:M \to M$.  (The involution $f$ is easy to  describe: send any
point $p \in M$ to the other point $q \in M$ that projects to that same point
under  $\pi$; if no such $q$ exists leave $p$ fixed.)  Using an invariant
version of Cerf theory we get $\s$ to intersect $f(\s)$ ``nicely'', and then
use $\s \cup f(\s)$ to construct a complex \ca\ with a list of useful
properties (Properties~\ref{prop-ca}).   \ca\ is used to construct an
invariant Heegaard surface for $M$ and bound its genus; the projection of this
surface gives the Heegaard surface for $N$, and estimating its genus we get a
linear upper bound for the genus of the quotient manifold in terms of $g(\s)$.

We now give the precise statements of our main results. 

\begin{rmkk}
\label{rmk:non-haken}{\rm
As is well-known, A Casson and C McA Gordon \cite{casson-gordon} proved that
if $M$ is an irreducible, non-Haken manifold then any minimal genus Heegaard
surface for $M$ is strongly irreducible.  Thus if $M$ is non-Haken we can
replace $g(\s)$ by $g(M)$ in all the statements below.} 
\end{rmkk}

\begin{thm}[Invariant Heegaard Surface]
\label{thm:first}
Let $M$ be \assu\  manifold \hgt\ admitting an \op\ involution $f$ and a
strongly irreducible \hhs\ $\s$. 

Then $M$ has an invariant \hhs\ of genus at most $\bddup.$
Moreover, each handlebody obtained by cutting $M$ open along this
surface is invariant.
\end{thm}

\begin{thm}[Genus of Double Covers]
  \label{thm:main}
Let $M$ be \assu\ manifold
\hgt\ admitting a strongly irreducible \hhs\ $\s$. Let $N$ be 
an orientable manifold and $\pi:M \to N$ a double cover.  Then we have:

$$g(N) \leq \bdddown.$$
\end{thm}

Using the invariant \hhs\ for $M$ constructed in Theorem~\ref{thm:first} we
obtain a bound on the complexity of the branch set. This bound is given in
terms of the bridge number of the branch set with respect to the Heegaard
surface for $N$ given in Theorem~\ref{thm:main}, \ie, the projection of the
invariant Heegaard surface for $M$.  The definition of bridge number with
respect to a Heegaard surface is given in Definition~\ref{dfn:bridge}  (for a
detailed discussion see, for example, \cite{kobayashi-rieck-growth} or
\cite{kobayashi-rieck-bridge}).  We prove: 

\begin{thm}
\label{thm:bridge}
Let $M$ be \assu\ manifold \hgt\ admitting a strongly irreducible \hhs\
$\s$. Let $N$ be  an orientable manifold and $\pi:M \to N$ be a double cover.
Denote the bridge index of the branch set with respect to the surface found in
Theorem~\ref{thm:main} by $b$.  Then we have: 

$$b \leq \bddupone.$$
\end{thm}

For proving Theorem \ref{thm:first} we study the intersection of strongly
irreducible Heegaard surfaces, that is, the intersection of $\s$ and its image
under the involution $f(\s)$.  However, our work can be applied for any two
strongly irreducible Heegaard 
surfaces $\s_1,\ \s_2 \subset M$.  We say that two embedded surfaces {\it
  intersect 
essentially} if every curve of $\s_1 \cap \s_2$ is essential in both
surfaces.  Rubinstein and M Scharlemann studied the intersection of strongly
irreducible Heegaard surfaces; we build on their work and prove:

\begin{thm}
\label{thm:m-cut-open}
Let $M$ be \assu\ manifold \hgt. Suppose that either $M$ admits
two strongly irreducible Heegaard surfaces $\s_1$ and $\s_2$ or a strongly
irreducible Heegaard surface $\s$ and an orientation preserving involution
$f$. Then we have:

\begin{enumerate}
  \item  $\s_1$ and $\s_2$ can be isotoped to intersect essentially
    and so that every component of $M$ cut open along $\s_1 \cup \s_2$ is a
    handlebody.
  \item  $\s$ can be isotoped so that $\s$ and $f(\s)$ intersect essentially
    and so that every component of $M$ cut open along $\s \cup f(\s)$ is a
    handlebody.
\end{enumerate}
\end{thm}

Theorem~\ref{thm:m-cut-open} follows quite easily from
Theorem~\ref{thm:spinal} (page~\pageref{thm:spinal}) which is
based on and improves results of Rubinstein and
Scharlemann, see Remark~\ref{rmk:improving-rs}.  We do not state this theorem
here to avoid terminology that had not yet been introduced.

Another tool used in the proof of Theorem~\ref{thm:first} is the creation of
an invariant complex $\ca\ \subset M$ fulfilling a list of properties
described in  \ref{prop-ca}.  Complexes fulfilling Properties~\ref{prop-ca}
are called {\it an-annular complexes}.  Properties~\ref{prop-ca} imply that
\ca\ has the following structure: it is constructed from a finite collection
of disjointly embedded tori (say $\{T_i\}_{i+1}^n$) bounding disjointly
embedded solid tori (say $\{V_i\}_{i+1}^n$) and a  collection of disjointly
embedded compact (but not closed) surfaces of {\it negative} Euler
characteristic\footnote{Getting the Euler characteristic of  this surface to
be negative (as opposed to non-positive)  is the main  challenge of the
construction and the reason for the name {\it an-annular}.}  with their
boundary on the tori $T_i$; all the boundary components form essential curves
on the tori. Properties~\ref{prop-ca} bound the Euler characteristic of \ca\
and state that $M$ cut open along \ca\ consists of handlebodies.  We refer
the reader to Section~\ref{sec:modify-ca} for a precise description of \ca\
and Properties~\ref{prop-ca}, and the statement and proof of
Theorem~\ref{thm:complex} where we prove the existence of \ca.

\begin{rmk}
Section~\ref{sec:modify-ca} is based on~\cite{rieck-proc} where Rieck proved
the existence of an-annular complexes in manifolds admitting two distinct
strongly irreducible Heegaard surfaces. 
\end{rmk}

Naturally, the solid tori  $\{V_i\}_{i+1}^n$ can be viewed as an equivariant
link in $M$.  Not every link in $M$ can be realized in this way and we ask
which links are (Question~\ref{que:link-exterior}).

This article is written in sections whose order, for the most
part, reveals the logic of the proof.  It is outlined in the next section.

\medskip

\noindent {\bf Acknowledgments:}   We are very grateful to Tsuyoshi Kobayashi and
Marc Lackenby for many helpful conversations and the anonymous referee for
her/his comments.  We thank Sean Bowman for the illustrations.

First named author: parts of this work was carried out while I was a JSPS post
doctoral fellow of Tsuyoshi Kobayashi in Nara Women's University, and while I
was visiting Akio Kawauchi in Osaka City University as a part of his 21st
Century COE Program ``Constitution of wide-angle mathematical basis focused on
knots''.   I am very grateful to both, and the math departments of Nara Women's
University and Osaka City University for their warm hospitality.

\section{outline}

\noindent {\bf Section \ref{s:background}, 
%page \pageref{s:background}:}    
pages \pageref{s:background}--\pageref{sec:solv}:}   
Background material, notation {\it etc.}

\smallskip

\noindent {\bf Section~\ref{sec:solv}, 
pages \pageref{sec:solv}--\pageref{sec:hhf}:}  
We give examples of higher order covers to demonstrate where our techniques
fail to generalize.  We also give examples that show the difficulty in finding
invariant reductions of various types (reducing sphere, weak reductions and
destabilizations).

\smallskip

\noindent {\bf Section~\ref{sec:hhf}, pages
  \pageref{sec:hhf}--\pageref{s:comp-free}:}  
We give a description of \hhf s (our version of sweepouts) and define the
Graphic.  Because of the invariance requirement the Graphic cannot be assumed
to be generic and this is rectified in Proposition~\ref{prop:crossing-graphic}
that shows that the behavior of the Graphic is essentially the same as the
behavior of generic graphics.  

\smallskip

\noindent {\bf Section~\ref{s:comp-free}, pages 
  \pageref{s:comp-free}--\pageref{s:spinal}:}  We isotope a
  strongly irreducible \hhs\ $\s \subset M$ to intersect its image under
  the involution in a compression free way, \ie, $\s$ and its image
  provide no compressions for each other, yet their intersection contains an
  essential curve.  In the end of this section we construct the generic
  interval, and isotopy of the Heegaard surface and its image that has the
  properties needed for Section~\ref{s:spinal}.  

\smallskip

\noindent {\bf Section~\ref{s:spinal}, pages
  \pageref{s:spinal}--\pageref{sec:m-cut-open}:} Using the generic interval we
  ensure $\s$ is chopped up completely by its image and a set of compressing
  disks for the image.  We also eliminate inessential \scc s of intersection
  between $\s$ and its image (that is, we isotope $\s$ to intersect its image
  essentially and spinally).  We also show that if  a manifold $M$ admits two
  strongly irreducible surfaces $\s_1$ and $\s_2$ (but not necessarily an
  involution), then $\s_1$ and $\s_2$ can be isotoped to intersect essentially
  and spinally. 

\smallskip

\noindent {\bf Section~\ref{sec:m-cut-open}, pages
  \pageref{sec:m-cut-open}--\pageref{sec:modify-ca}:}  Proof of   Theorem
  \ref{thm:m-cut-open}.

\smallskip

\noindent {\bf Section~\ref{sec:modify-ca}, pages
  \pageref{sec:modify-ca}--\pageref{sec:get-s}:} We consider  $\s$
  union its  image as a complex and modify it to get rid of undesired annuli,
  proving (Theorem~\ref{thm:complex}) existence of the complex \ca\ fulfilling
  Properties~\ref{prop-ca}.

\smallskip

\noindent {\bf Section~\ref{sec:get-s}, pages
  \pageref{sec:get-s}--\pageref{sec:hhs-for-quotient}:} Using this complex we
  create an   invariant \hhs\ for $M$ and estimate its genus, thus proving
  Theorem \ref{thm:first}.

\smallskip

\noindent {\bf Section~\ref{sec:hhs-for-quotient}, 
%pages \pageref{sec:hhs-for-quotient}--\pageref{sec:bridge}:} 
page \pageref{sec:hhs-for-quotient}:} 
Using this \hhs\ and the Equivariant Disk Theorem we get a \hhs\ for the
quotient thus proving Theorem \ref{thm:main}.

\smallskip

\noindent {\bf Section~\ref{sec:bridge}, pages
  \pageref{sec:bridge}--\pageref{theEnd}:} Using the bounded genus invariant
  surface found in Section \ref{sec:get-s} we bound the complexity of the
  branch set in terms of the genus of $M$.

\section{Background}
\label{s:background}

We work in the smooth and orientable category.  By {\it manifold} we mean a
3-dimensional compact manifold without boundary.  We follow standard notation
for 3-manifolds: $\mbox{int}X$ is the interior of $X$, $\mbox{cl}X$ is the
closure of $X$, $\del X$ is the boundary of $X$ {\it etc.} 
See \cite{hempel} or \cite{jaco} for basic definitions.  We refer the reader to
\cite{scharlemann-survey} for a detailed discussion about Heegaard
splittings.  We assume that our manifold
$M$ is not a \sfs. We note that for \sfs s results far more
refined than ours are known, \eg\ for $S^3$ the positive solution
of the Smith Conjecture \cite{smith}, Hodgson and Rubinstein's work
about lens spaces \cite{hodgson-rubinstein},  Boileau and Otal's
work about small \sfs s, \cite{MR92i:57014}, and 
Scott's work about Haken \sfs s \cite{MR87c:57012}.

We further assume that our manifold contains a strongly
irreducible \hhs, \ie, $M$ has a \hhs\ for which any two
compressing disks on opposite sides intersect. By Haken
\cite{MR36:7118} our manifold is irreducible, \ie, every 2-sphere
embedded in $M$ bounds a ball. This condition is not vacuous: Casson and
Gordon's \cite{casson-gordon} seminal work show that for irreducible 
non-Haken manifolds every
minimal genus (indeed, any irreducible) \hhs\ is strongly
irreducible. (Non-Haken manifolds are not the only manifolds that
contain strongly irreducible Heegaard splittings; see T Kobayashi and Rieck
\cite{kobayashi-rieck-reductions} for manifolds admitting both weakly
reducible and strongly irreducible minimal genus Heegaard splittings.)  We
note however that constraints are imposed on the cover and not on the
manifold being covered, where no additional constraints apply.

Suppose that $p:M \to N$ is a double cover.  Since all double
covers are regular (including branched double covers), there
exists $f:M \to M$ an involution on $M$, so that $p$ is given by
the natural projection $M \to M/(f) \cong N$.  A subset $S \subset M$ is called
{\it invariant} if $f(S)=S$.  $S \subset  M$ is called {\it equivariant} if
$S$ is either invariant or disjoint from its own image, \ie,
either $f(S)=S$ or $f(S) \cap S = \emptyset$.  We use the notation
$N(S)$ to mean a normal neighborhood.  When discussing an
invariant (resp.  equivariant) subset of $M$, we use $N(S)$ to
denote an invariant (resp. equivariant) normal neighborhood.

Since all manifolds are assumed to be orientable, $f$ is
orientation preserving.

The first half of this paper deals with the intersection of embedded
surfaces.  We follow the terminology used in \cite{rs1:1996}.  In particular,
a tangency between two embedded surfaces can have one of two forms: a {\it
  center} (modeled on the intersection of $z=0$ and $z=x^2+y^2$) or saddle
(modeled on the intersection of $z=0$ and $z=x^2-y^2$ at the origin).

\section{Examples}
\label{sec:solv}

Now an example.  We consider Solv manifolds (definition below, see
also \cite{scott:1983}) since they cover each other generously and
in \cite{cooper-scharlemann} Cooper and Scharlemann gave a
complete classification of their \hhs s.

\begin{dfn}
\label{dfn:solv}
A 3-manifold is called {\it Solv} if it is a torus bundle over $S^1$ with
Anosov monodromy, \ie, the monodromy has infinite order and no power 
of it has fixed point in $\pi_1(T^2)$. 
\end{dfn}

Given a Solv manifold $M$ (say with monodromy $\phi$) and a positive integer
$n$, the Solv manifold with monodromy $\phi^n$ (denoted $M_n$) is an n-fold
cover of $M$.   This cover is as nice as one could hope for: cyclic (in
particular regular) and unbranched. As it is our goal to get invariant \hhs s
(Theorem~\ref{thm:first}), it is interesting to consider a minimal genus \hhs\
for $M_n$, say $\s$.  $\s$ and its image under a generator of the action of
$\mathbb{Z} /(n)$ are \hhs s of genus 2 or 3.   By picking $M_n$ correctly,
\cite{cooper-scharlemann} show that $\s$ and its image under the generator of
the action are isotopic.  But for a cyclic group action invariance under a
generator implies invariance under the entire group; may we conclude that the
surface is invariant?

No.  For if it were, for larger and larger values of $n$ we would get
surfaces of genus 2 or 3 that are invariant under the free action
of a cyclic group of arbitrarily high order.  This of course
cannot be, since the quotient surface would have fractional Euler
characteristic (so a surface invariant under a group action of high order must
have high genus).  $\s$ ``equals'' its image in the sense of
``up-to-isotopy'', but this isotopy cannot be realized invariantly and
therefore it does not provide us with a surface that is truly invariant under
the group action.  Moreover, it is an easy exercise to get $\s$ to be disjoint
from its image.  Yet during the isotopy that takes $\s$ to its image the two
are no longer each other's images.   This phenomenon occurs in a very simple
setting as well: consider $S^1$ double covering itself.  The preimage of a
point is two points which are isotopic to each other but the isotopy cannot
be  realized invariantly.  
This example demonstrates that generalizing this work for higher
order coverings will not be a straight forward task but will
require a new ingredient, perhaps the degree of the cover.

We use Solv manifolds to provide one more example.  Suppose $M$ is a Solv
manifold of genus two, $\s \subset M$ a minimal genus Heegaard surface,
and $M_n$ its n-fold cyclic cover (as above).  Since $M$ is
irreducible $\s$ is strongly
irreducible.   Lifting $\s$ to $M_n$ we get a Heegaard surface of genus $n +
1$, say $\s_n$.  By \cite{cooper-scharlemann} the only irreducible Heegaard
surfaces for $M_n$ are minimal genus Heegaard surfaces.  Therefore $\s_n$
destabilizes $n-1$ or $n-2$ times, which gives many distinct collections of
reducing, weakly reducing and destabilizing disks for $\s_n$.  However, strong
irreducibility of $\s$ implies that no such reducing set can be made
equivariant under the cyclic group action.  Therefore, either on at least one
side the disks are not equivariant, or the disks are equivariant on both sides
but in the projection of the disks to $M$ every pair of disks on opposite
sides of $\s$ intersect at least twice (note that the image of disjoint
curves on $\s_n$ may intersect more than once).

\section{\hh\ Functions and the Graphic}
\label{sec:hhf}

In this section we introduce the basic set up, beginning with the
following definition that formalizes the basic tool we use for
studying \hhs s.  It is equivalent to the notion of {\it
sweepout}, as defined in \cite{rs1:1996}.  (Much of the material
in this section is not new but is included here for our work in the
equivariant setting.)

\begin{dfn}
\label{d:heegaard-function} Let $M$ be a manifold.  A smooth function
$h:M \to [-\infty,\infty]=\{-\infty\}\cup \bbb R \cup \{\infty\}$ is called a
{\it \hhf\ } if the following hold:

   \begin{enumerate}
    \item $h^{-1}(-\infty)$ and $h^{-1}(\infty)$ are graphs;
    \item $h|_{h^{-1}(\bbb R)}$ has no critical points.
   \end{enumerate}
\end{dfn}

(2) implies that for any $t \in \mathbb{R}$ $h^{-1}(t)$ is a smooth surface 
and its genus is independent of $t$.  In fact,
any two such surfaces (say $h^{-1}(s)$ and $h^{-1}(t)$, with
$s<t$) are parallel and the region defining the parallelism is
given by $h^{-1}([s,t])$.  

\begin{dfns}
\label{dfn:hhs}
  \begin{enumerate}
  \item A surface $\s \subset M$ is called a {\it \hhs} if it is of
    the form $h^{-1}(0)$ for some \hhf\ $h$.
  \item A {\it spine} for a \hhs\ $\s$ is a (disconnected) embedded graph of
  the form $h^{-1}(-\infty) \sqcup h^{-1}(\infty)$.
  \end{enumerate}
\end{dfns}

We will often start with a \hhs, and then consider a \hhf\ that gave rise
to it, \ie, we shall start with $\s$ and consider $h$ as in
Definition \ref{dfn:hhs}. This function will be called a
``corresponding \hhf''.  It is by no means unique, nor is the spine.

Let $M$ be a manifold and $f:M \to M$ an involution.
Let $\s$ be a \hhs\ for $M$, and $h$ a corresponding \hhf.  When
studying $\s$ and $f(\s)$ the \hhf\ corresponding to $f(\s)$
we will use is $h \circ f$ (note that $f=f^{-1}$).

We will use a Cerf theoretic argument, which requires the
construction of the Graphic.  The Graphic is based on a
2-parameter family of surfaces, \ie, the assignment of two
surfaces for every point in the parameter square $\{(s,t):s,t\in
[-\infty,\infty]\}$, denoted $(s,t) \mapsto (F_1(s,t),F_2(s,t))$.
The Graphic itself is the subset of points corresponding to
surfaces that do not intersect transversely. See \cite{rs1:1996}
or \cite{rieck:1997.2} for further details about the Graphic, or
Cerf's original work \cite{cerf}.  (We give a more detailed
description of the graphic below.)

In our case, given $h$ a \hhf\ for $M$ and $f$ an involution on
$M$, we start with the assignment: $(s,t) \mapsto (\s_s,\gst)$
where $\s_s=h^{-1}(s)$ and $f(\s_t)=f(h^{-1}(t))$.  Note that on
the diagonal $\{ s=t \}$ the involution exchanges the two surfaces.
This assignment is not necessarily generic and therefore no
niceness properties of the graphic can be assumed (not even
one-dimensionality).  To that end, we modify the surfaces. First,
and most importantly, we modify the surfaces along the diagonal,
as in \cite{hodgson-rubinstein}.  Via perturbation, we impose the
following two conditions on $h$: the spines $h^{-1}(\pm\infty)$
are disjoint from their images, and zero must be a regular value
of $h - h \circ f$.\footnote{This condition is quite natural: we 
are interested in the intersection of the surface $h^{-1}(t)$ with 
its image; therefore we are forced to look at points where $h$ and $h \circ 
f$ have the same value.} 
Since both conditions are generic after imposing the first
condition on the spines we can impose the second without ruining
the first.  In \cite{hodgson-rubinstein} it was shown that these
conditions imply the following:

\begin{conds}
\label{conds:generic}
\begin{enumerate}
\item For all but finitely many values of $t$ the intersection of
  $\s_t$ and $\gst$ is transverse.  Points that correspond to non-transverse
  intersection are called {\it critical}. At a critical point exactly one of
  the following holds: 
\item $\s_t$ and $\gst$ intersect in a single non-degenerate critical
  point that is fixed by $f$ (then $t$ is called a {\it simple} critical
  point).
\item $\s_t$ and $\gst$ intersect at a pair of non-degenerate critical points 
  that are exchanged by $f$ (then $t$ is called a {\it double} critical
  point). 
\end{enumerate}
\end{conds}

Before modifying the 2-parameter family off the diagonal to get
the surfaces to be as generic as we can, let us explain what we
mean by genericity.  It is a local property, \ie, given a point
$(s_0,t_0)$ it only depends on surfaces for $(s,t)$ close to
$(s_0,t_0)$.  In particular we can impose it on an open set, which
we shall do (the complement to the diagonal, to be precise).  The
condition is this: for a dense open set of points $(s,t)$ the
surfaces intersect transversely.  The points where the surfaces
intersect non transversely fall into three categories.

\smallskip

\noindent {\bf Edges:} 1-dimensional sets in the parameter square, with
  finitely many components each homeomorphic to an interval.
  The points of the edges are those that correspond to pairs of surfaces
  having exactly one critical point.

\smallskip

\noindent {\bf Vertices of valence four:} finitely many
  points in the parameter square that correspond to pairs of surfaces with
  exactly two critical points.  Each valence four vertex is the endpoint of
  exactly  four edges, more precisely, two pairs of edges where each pair
  corresponds to a tangency between the two surfaces.  We can also
  consider each pair of edges as one long edge, pasting them together at the
  vertex.  Then the vertex is the point where the two edges cross each other
  transversely.

\smallskip

\noindent {\bf Death-birth vertices:} finitely many vertices of valence two.
  As they play no role
  whatsoever in this work, so we do not describe them here.

So a generic Graphic forms a finite graph embedded in the parameter square.
We now prepare the Graphic: starting with $(s,t)
\mapsto (\s_s,\gst)$, on one side of the diagonal (say $s>t$) we perturb the
surfaces to be generic.  We may do so without changing the diagonal: we take
any generic perturbation of the surfaces at $s>t$, say given by $(s,t) \mapsto
(F_{1,r}(s,t),F_{2,r}(s,t))$ so that at $r=0$ we have our original assignment,
and we pick a perturbation given by $(s,t) \mapsto
(F_{1,r'}(s,t),F_{2,r'}(s,t))$, with $r'$ a function of $(s,t)$ that limits on
zero as $(s,t)$ approaches the diagonal.  If all the above is done
generically, we have an assignment that is generic at $s>t$, fulfills
Conditions~\ref{conds:generic} on the diagonal, and is continuous on $s \geq
t$.  For $s<t$, we set $F_1(s,t) = F_2(t,s)$ and $F_2(s,t)=F_1(t,s)$.  (In
other words, we perturb the surfaces in the domain $s<t$ in the exact same way
we did in the domain $s>t$.)  Note that $f$ exchanges $F_1$ and $F_2$:
$f(F_1(s,t)) = F_2(t,s)$ and $f(F_2(s,t)) = F_1(t,s)$ .  Hence $f$ induces the
involution $(s,t) \mapsto (t,s)$ on the parameter square and the Graphic is
invariant under this involution (in general, this forces double critical
points on the diagonal).

\begin{rmkk}
\label{rmkk:notation} {\rm We work mostly on the diagonal, where
the surfaces are parameterized by a single parameter $t$,
explicitly: $\s_t=F_1(t,t)$ and $\gst=F_2(t,t).$  For simplicity,
while considering points on the diagonal we use the notation
$\s_t$ and $\gst$.  }
\end{rmkk}

We conclude this section with the following proposition, which discusses the
behavior of the Graphic near the diagonal.  When the Graphic intersects the
diagonal, this proposition is needed due to lack of transversality.  Here and
throughout this work we move freely between a point $(s,t)$ of the parameter
square and the corresponding surfaces $F_1(s,t)$ and $F_2(s,t)$, and between
edges on the Graphic and the corresponding tangencies of $F_1(s,t)$ and
$F_2(s,t)$.

\begin{pro}
\label{prop:crossing-graphic}
Let $(t_0,t_0)$ be a point on the
diagonal that corresponds to surfaces with two critical points,
and let $S_1$ and $S_2$ be the two curves of the graphic through
it, each corresponding to one of the critical points.
Then locally about $(t_0,t_0)$ one of the following holds:

\begin{enumerate}
\item $S_1$ is on one side of the diagonal (except for
$(t_0,t_0)$)
  and its image $S_2$ on the other;
\item $S_1$ and $S_2$ cross each other.
\end{enumerate}
\end{pro}

\begin{proof}
By Conditions \ref{conds:generic}, locally the curve $S_1$ has only one point
on the diagonal.  Assume (1) does not occur.  Therefore, one of the two curves
(say  $S_1$) crosses the diagonal.  If $S_1=S_2$ near $(t_0,t_0)$ then the two
critical points are in fact the same (since off the diagonal the graphic is
generic), contrary to our assumption. The proposition follows from the fact
that $S_2$ is the image of $S_1$ under the involution $(s,t) \to (t,s).$
\end{proof}

\begin{rmks}
  \begin{enumerate}
  \item In the first case, while traveling along the diagonal, it is as
    if we did not encounter a critical point at all, as the diagonal
    is only tangent to the two but never traverses them.
    We may ignore such points throughout this work.
  \item In the second case the intersection behaves as if the graphic
    is generic.
  \end{enumerate}
\end{rmks}

\section{Compression Free Intersection}
\label{s:comp-free}

Following Rubinstein and Scharlemann \cite{rs1:1996} we define:

\begin{dfns} 
\label{d:comp-free} Let $F_1$ and $F_2$ be surfaces embedded in a
3-manifold intersecting transversely.

\begin{enumerate}
\item A curve of $F_1 \cap F_2$ is called {\it essential}
(resp. {\it inessential}) if it is essential (resp. inessential)
in both $F_1$ and $F_2$.  A  curve of $F_1 \cap F_2$ which is essential on
one surface and inessential on the other is called a {\it compression}.

\item  The intersection of $F_1$ and $F_2$ is called {\it
compression free} if no curve of $F_1 \cap F_2$ is a compression.
\end{enumerate}
\end{dfns}

Recall (Remark~\ref{rmk:non-haken}) that if $M$ is non-Haken then any minimal
genus \hhs\ is strongly irreducible; hence the theorem below is not vacuous:

\begin{thm}
\label{thm:compr-free} Let $M$ be \assu\ manifold \hgt\ with an
\op\ involution $f:M \to M$. Let $\s$ be a strongly irreducible
\hhs\ for $M$, and $h:M \to \mathbb R$ a corresponding \hhf. Then
there exists an interval $(a,b)\subset \bbb R$ so that the following
conditions hold:

\begin{enumerate}
\item For any regular point $t \in (a,b)$ the intersection of
$\s_t=h^{-1}(t)$ with $\gst$ is compression free, yet contains an
essential curve.
\item $a$ is critical, and arbitrarily close to $(a,a)$ there are points
  corresponding to transverse intersections that are not compression free or
  are all inessential\footnote{These points may be off the diagonal.}
  (similarly for $b$). 
\item For any $t \in (a,b)$, $(t,t)$ has a
  neighborhood $U$ so that every regular point in $U$
  corresponds to a compression free intersection containing an essential
  curve.  
\end{enumerate}
\end{thm}

\begin{proof}
Color the handlebody $h^{-1}([-\infty,s])$ purple and $h^{-1}([s,\infty])$
yellow. 

For convenience we present the diagonal as an interval by
identifying any point $(t,t)$ on the diagonal with $t \in \mathbb
R$.  Subdivide the interval $[-\infty,\infty]$ into {\it layers}
separated by critical points (recall their definition in
Conditions~\ref{conds:generic}).  A component of the parameter square cut open
along the Graphic is called a {\it region}.  Note that every layer is
contained in exactly one region.  Since the intersection pattern
between $\s_t$ and $\gst$ is independent of the choice of point
within a region (resp. a layer), we call it {\it the intersection of the
  region} (resp. the layer). We label regions and layers, exactly as in
\cite{rs1:1996}.  By construction $F_1(s,t)$ is a small perturbation of
$\s_s$ and therefore the two handlebodies complementary to
$F_1(s,t)$ inherit yellow and purple coloring.

\begin{dfns}
\label{dfn:labels}

Let $R$ be a region of the Graphic and $(s,t) \in R$, with
corresponding surfaces $F_1(s,t)$ and $F_2(s,t)$.

     \begin{enumerate}
     \item $R$ is labeled \bfP\ if there exists a disk
     $D_{\bfP}\subset F_2(s,t)$ and the following conditions holds:
               \begin{enumerate}
               \item The boundary of $D_{\bfP}$ is an essential curve of
               $F_1(s,t)$. 
               \item Near its boundary $D_{\bfP}$ is purple.
               \item $\mbox{int}D_{\bfP} \cap F_1(s,t)$ does not contain
                 essential curves of $F_1(s,t)$.
               \end{enumerate}
     \item $R$ is labeled   \bfp\ whenever  the following
     conditions holds:
               \begin{enumerate}
               \item The intersection of $F_2(s,t)$ and the yellow handlebody
               contains an essential curve of $F_2(s,t)$.
               \item Every curve of $F_1(s,t) \cap F_2(s,t)$ is
               inessential.
               \end{enumerate}
     \item The labels \bfY\ and \bfy\ are defined similarly.
     \item A layer is labeled by the same labels as a region that
     contains it.
     \end{enumerate}
\end{dfns}

\begin{rmk}
Rubinstein and Scharlemann define their labels in \cite[pages
1009--1010]{rs1:1996}.  To see that their labels are indeed the
same as ours, denote $F_1(s,t)$ by $P$, $F_2(s,t)$ by $Q$, the purple
handlebody by $A$ and the yellow handlebody by $B$.  Then the
labels \bfP, \bfp, \bfY, and \bfy,  correspond to the labels $A$,
$a$, $B$, and $b$ in \cite{rs1:1996} (in the same order).  The
label $X$, $x$, $Y$ and $y$ appearing in \cite{rs1:1996} are not
needed here (essentially, because the surface $\gst$ is the image of
$\s_t$ and has the same intersection properties).
\end{rmk}

From \cite[Sections~4 and 5]{rs1:1996} we have the proposition below for
regions of the Graphic, that is applicable directly for simple critical points
of the diagonal, as layers separated by simple points are contained in regions
that share an edge.  The goal of this section  can be described as extending
this proposition to double critical points (we remark that in general this
cannot be done, and we will need to use the involution).

\begin{pro}
\label{pro:adjacent-regions} Every region of the Graphic has at
most one label. A region labeled \bfp\ or \bfP\ cannot share
an edge with a region labeled \bfy\ or \bfY.  Therefore every layer has at
most one label and a layer labeled \bfp\
or \bfP\ cannot be separated by a simple critical point from a layer labeled
\bfy\ or \bfY. 
\end{pro}

In Lemmas~\ref{lem:lowercase} and \ref{lem:bfpXbfy} we study
lowercase labels:

\begin{lem}
\label{lem:lowercase}
A region $R$ is labeled with a lowercase label if and only if
the corresponding surfaces intersect in inessential curves only.
In that case, except for punctures the surface $F_2(s,t)$ has one
color, purple if the label is \bfy\ and yellow if the label is
\bfp.
\end{lem}

\begin{proof}
This information is contained in \cite{rs1:1996} so we paraphrase
it here.  If a region has a lowercase label then by definition the intersection
consists of entirely inessential curves.  Conversely, if
the intersection consists of inessential curves only, except
perhaps for punctures $F_2(s,t)$ is colored in one color, yellow
or purple (resp.).  Any essential curve of $F_2(s,t)$ can be
isotoped off the punctures, showing the label is \bfp\ or \bfy\
(resp.).
\end{proof}

\begin{lem}
\label{lem:bfpXbfy} There does not exist a critical level $t_0$ separating a 
layer labeled \bfy\ from a layer labeled \bfp.
\end{lem}

\begin{proof}
Suppose for contradiction there exists a critical level $t_0$ that separates a
layer $l_y$ labeled \bfy\ from a layer $l_p$ labeled \bfp.  By
Proposition~\ref{pro:adjacent-regions}, $t_0$ is a double critical point.  It
is easy to see that \bfp\ does not change when crossing centers, hence $t_0$
is a double saddle.  By Lemma~\ref{lem:lowercase}, for $t_y \in l_y$ exactly
one component of $f(\s_{t_y}) \cap h^{-1}([-\infty,t_y])$ (say $F$) contains all of
$f(\s_{t_y})$ except, perhaps, for punctures.  After the double saddle, no
essential curve of $f(\s_{t_p})$ is contained in $h^{-1}([-\infty,t_p])$.  Since
crossing each saddle changes $f(\s_{t_y}) \cap h^{-1}([t_y,\infty])$ by adding or
removing a single 1-handle, we see that all of $F$ (except perhaps for
punctures) was moved out of $h^{-1}([-\infty,t_y])$ and into $h^{-1}([t_p,\infty])$ by  two
1-handles; hence, $F$ is a punctured torus (resp. pair of pants) and $\s_t$
is a torus (resp. sphere), contradicting the assumptions of
Theorem~\ref{thm:compr-free}.
\end{proof}

In Lemma~\ref{lem:PnearY} we study uppercase labels:

\begin{lem}
\label{lem:PnearY} There does not exit a critical level $t_0$
separating a layer labeled \bfY\ from a layer labeled \bfP.
\end{lem}

\begin{proof}
Assume for contradiction $t_0$ is a critical level separating a layer labeled
\bfP\ (say $l_p$) from a layer labeled \bfY\ (say $l_y)$.  For convenience of
presentation we assume that $l_p$ is below $t_0$ and $l_y$ is above $t_0$
(the other case is treated by taking $\epsilon$ below to be a small negative
number).  By
Proposition~\ref{pro:adjacent-regions} $t_0$ is a double critical point. As
above, it is easy to see that uppercase labels do not change when crossing
centers; hence we may assume that $t_0$ corresponds to two saddles (say $S_1$
and $S_2$).  From Definition~\ref{dfn:labels} we see that for  $t \in l_p$,
there exists a compressing disk (say $D_{\bfP} 
\subset \gst$) giving rise to the label \bfP. Similarly, for $t \in l_y$ there
exists a compressing disk $D_{\bfY}$ giving the label \bfY.  We consider two
cases:

\medskip

\noindent {\bf Case 1.  One of the saddles does not destroy one of
the disks ($D_{\bfY}$ or $D_{\bfP}$):}  Leaving the diagonal and
crossing the edge we hypothesized not to destroy one of the disks,
we see that a region labeled $\bfY$ is adjacent to a region labeled $\bfP$ 
along an edge of
the graphic, contradicting Proposition~\ref{pro:adjacent-regions}.

\medskip

\noindent {\bf Case 2.  Each saddle $S_1$ and $S_2$ destroys both
$D_{\bfY}$ and $D_{\bfP}$:} Let $\gamma = \del D_{\bfP}$.  As we
approach $t_0$ the curve $\gamma$ limits on both saddles, or we
would be in Case~1. Since the saddles are involutes of each other,
$f(\gamma)$ must limit on both saddles as well.  There are two
subcases:

\medskip

\noindent {\bf Subcase 2.a:} $\gamma = f(\gamma)$.  Since
$\gamma=\del D_{\bfP}$ and $D_{\bfP} \subset \gst$, we see that
$\del f(D_{\bfP}) = f(\del D_{\bfP}) = f(\gamma) = \gamma$.  But
$f(D_{\bfP})\subset f(\gst) = \s_t$.  Therefore $\gamma$ is
inessential in $\s_t$, contradicting Definition~\ref{dfn:labels}. 

\medskip

\noindent {\bf Subcase 2.b:} $\gamma \neq f(\gamma)$. Since
$f(\gamma)$ must  limit on both saddles as well, 
$\gamma$ limits on  each saddle once only.

Let $\epsilon>0$ be a small number.
Let $v$ be a non-vanishing vector field along $\gamma$  that is everywhere
transverse to $\gamma$ and is tangent to $\s_{t_0-\epsilon}$.  We may assume
$v$ points towards $S_1$ as $\epsilon$ tends to zero 
(else we reverse it).  The curve $\gamma$ is called {\it untwisted} if $v$
points towards $S_2$ as $\epsilon$ tends to zero, {\it twisted} otherwise
(\ie, if $v$ points away from $S_2$ as $\epsilon$ tends to zero); note that
this is independent of choice of $v$. Likewise, $f(\gamma)$ may 
be twisted or untwisted; when examining twistedness of $f(\gamma)$ we regard
it as a curve on $\s_{t_0-\epsilon}$, not on $f(\s_{t_0-\epsilon})$. We show
that $\gamma$ is twisted if and only if $f(\gamma)$ is: suppose $\gamma$ is
untwisted.  Near $S_1$ and $S_2$ we can view $\s_{t_0-\epsilon}$ as a flat
disk in $z=-\epsilon$, and $f(\s_{t_0-\epsilon})$ as a small piece of
$z=x^2-y^2$ forming a little arch. 
We denote the handlebodies given by $M$ cut open along
$f(\s_{t_0-\epsilon})$ by $H_1$ and $H_2$.  By construction $v$ is tranverse to
$f(\s_{t_0-\epsilon})$ and after renaming $H_1$ and $H_2$ if necessary we may
assume that $v$ points out of $H_1$.  Let $u$ be
the vector field along $f(\gamma)$ that points out of $H_1$ and is tangent to
$\s_{t_0-\epsilon}$.  By assumption, near $S_1$ and $S_2$ $v$ points towards
the saddles or ``into'' the arches; hence $H_1$ is above both arches and $u$
points into the arches as well.  Thus $u$ points towards both
saddles, showing that $f(\gamma)$ is untwisted.  Thus $\gamma$ is twisted
implies that $f(\gamma)$ is twisted and similarly we see that $f(\gamma)$ is
twisted implies that $\gamma$ is twisted, as required.

After crossing $S_1$ $\gamma$ and $f(\gamma)$ form a single curve (say
$\beta$) and after crossing $S_2$ this curve breaks up into two involute
curves, say $\alpha$ and $f(\alpha)$.  As we approach $t_0$ from above  the
boundary of $D_{\bfY}$ limits on both saddles.  Thus the boundary of
$D_{\bfY}$ is $\alpha$ or $f(\alpha)$.

Assume first both $\gamma$ and $f(\gamma)$ are  untwisted and let $v$
(resp. $u$) be a vector field along $\gamma$ (resp. $f(\gamma)$) pointing
towards both saddles.  Pushing $\gamma$ (resp. $f(\gamma)$) slightly along $v$
(resp. $u$)  we obtain a curve $\gamma'$ (resp. $\gamma''$).  Performing both
saddle crossings on $\gamma'$ and $\gamma''$ we get two curves that are
isotopic to $\alpha$ and $f(\alpha)$ and are disjoint from $\gamma$ and
$f(\gamma)$.  See Figure~\ref{fig:reductions}.
      \begin{figure}
      \centerline{  \includegraphics[width=4in]{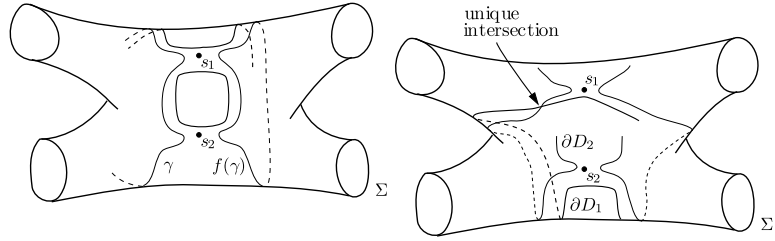}}
      \caption{The reductions: untwisted (left) and twisted (right)}
      \label{fig:reductions}
      \end{figure}
Although the disks $D_{\bfP}$ and $D_{\bfY}$ may intersect the Heegaard
surface in their interior, by definition of uppercase labels any such curve of
intersection is trivial in the Heegaard surface and (as noted in
\cite{rs1:1996})  the boundary of $D_{\bfP}$ bounds a purple meridian disk and
the boundary of $D_{\bfY}$ bounds a yellow meridian disk.  We conclude that
$\s_{t_0}$ weakly reduces, contradicting our  assumption.

Next assume that both $\gamma$ and $f(\gamma)$ are twisted.  Let $v$ and $u$
be  vector fields pointing towards $S_1$ and away from $S_2$.  Again, we push
$\gamma$ along $v$ and $f(\gamma)$ along $u$ obtaining $\gamma'$ and
$\gamma''$.  Performing the saddle operation $S_1$ on $\gamma'$ and
$\gamma''$, we obtain the curve $\beta$ and see that $\beta$ is disjoint
from $\gamma$ and $f(\gamma)$.   Performing the  saddle operation $S_2$ on
$\beta$ we obtain two curves (say $\alpha'$ and $\alpha''$).  Since $\gamma$
and $f(\gamma)$ separate $\beta$ from $S_2$, the curves obtained are not
disjoint from $\gamma$.  However, it is easy to see directly that $|\alpha'
\cap \gamma|  = 1$ and $|\alpha'' \cap \gamma| = 1$.  As above, $\gamma$
bounds a purple meridian disk and either $\alpha'$ or $\alpha''$ is  isotopic
to $\del D_{\bfY}$ and hence bounds a yellow meridian disk.  We conclude that
$\s_t$ destabilizes, contradicting our assumptions.
\end{proof}

In Lemma~\ref{lem:P-near-y} we study mixed labels:

\begin{lem}
\label{lem:P-near-y}
There does not exit a critical level $t_0$ separating a layer labeled \bfy\
from a layer labeled \bfP\ (and similarly for \bfY\  and \bfp).
\end{lem}

\begin{proof}
Assume for contradiction $t_0$ is a critical level separating a layer labeled
\bfP\ (say $l_p$) from a layer labeled \bfy\ (say $l_y)$.  For convenience we
assume $l_p$ is below $l_y$ (the other case is treated by taking $\epsilon$
below to be a small negative number).  For $t \in l_p$ there exists a
compressing disk  (say $D_{\bfP}$) that gives rise to the label \bfP.  As in
Lemma~\ref{lem:PnearY} we may assume that $t_0$ corresponds to two involute
saddles, say $S_1$ and $S_2$. We consider two cases:

\medskip

\noindent {\bf Case 1:  One of the saddles does not destroy
$D_{\bfP}$.}  This is identical to Lemma~\ref{lem:PnearY}~(1).   We may assume
from now on this is not the case.

\medskip

\noindent {\bf Case 2:  Both $S_1$ and $S_2$ destroys $D_{\bfP}$.}
Let $\epsilon>0$ be sufficiently small.  For $t_0 - \epsilon$, let $\gamma =
\del D_{\bfP}$.  As $\epsilon$ approaches zero the curve $\gamma$ limits on
both saddles, or we would be in case~1. Since the saddles are involutes of
each other, $f(\gamma)$ must limit on both saddles as well.  There are two
subcases:

\medskip

\noindent {\bf Subcase 2.a:} $\gamma = f(\gamma)$.   This is
identical to Lemma~\ref{lem:PnearY}~(2.a).  We may assume from now on this is
not the case.

\medskip

\noindent {\bf Subcase 2.b:} $\gamma \neq f(\gamma)$.  On $\s_{\tpe}$ all
curves of intersection are inessential.  Denote the region containing $l_y$ by
$R_y$ (similarly $R_p$) and the region we get to after crossing $S_1$ out of
$R_y$ by $R$.  Since $R$ shares an edge $R_y$ and another edge with $R_p$, by
Proposition~\ref{pro:adjacent-regions} $R$ is unlabeled.  Fix $(s,t) \in R$.
If $F_1(s,t) \cap F_2(s,t)$ consists entirely 
of inessential \scc s then by Lemma~\ref{lem:lowercase} $R$ has a lowercase
label, contradiction.   Suppose that $F_1(s,t) \cap F_2(s,t)$ contains curves
that are essential in $F_1(s,t)$ (the other case is symmetric).  Thus we see
that crossing $S_1$ a single inessential   curve of intersection becomes two
\scc s that are both essential in $F_1(s,t)$ say $\alpha'$ and $\alpha''$.
Note that $\alpha'$ is parallel to $\alpha''$ in $F_1(s,t)$ and all other
curves of $F_1(s,t) \cap F_2(s,t)$ are inessential in $F_1(s,t)$.

By symmetry of $\alpha'$ and $\alpha''$, there are four possibilities when
crossing $S_2$ out of $R$ into $R_p$:

\begin{enumerate}
\item $S_2$ does not involve $\alpha'$ or $\alpha''$.
\item $S_2$ connects $\alpha'$ to a \scc\ that is inessential in $F_1(s,t)$.
\item $S_2$ connects $\alpha'$ to itself.
\item $S_2$ connects $\alpha'$ to $\alpha''$
\end{enumerate}

We conclude the proof of Lemma~\ref{lem:P-near-y} by reducing~(1)--(4) above
to previous subcases: in~(1), no curve of $\s_{t_0 - \epsilon} \cap f(\s_{t_0
- \epsilon})$ involves both saddles hence this is in fact Case~1 above,
contradiction.  In~(2), let $\beta$ be the  curve obtained from $\alpha'$ and
an inessential \scc\ after crossing $S_2$.  Then $\beta$ is the unique curve
of  $\s_{t_0 - \epsilon} \cap f(\s_{t_0 - \epsilon})$ involving both saddles.
In~(3) the curve $\alpha'$  splits into two curves and exactly one of the two
involves both saddles.  In~(4) the curves $\alpha'$ and $\alpha''$ become a
single curve.  Thus, in~(2)--(4) there is a unique  curve of
$\s_{t_0 - \epsilon} \cap f(\s_{t_0 - \epsilon})$ that involves both saddles;
since both $\gamma$ and $f(\gamma)$ involve both saddles, we conclude that
$\gamma = f(\gamma)$ and we are in fact in Subcase~2.a.  With this
contradiction we conclude the proof of Lemma~\ref{lem:P-near-y}.
\end{proof}

Combining Proposition~\ref{pro:adjacent-regions} and Lemmas~\ref{lem:bfpXbfy},
\ref{lem:PnearY} and \ref{lem:P-near-y} we get:

\begin{pro}
\label{pro:adjacent-layers} A layer labeled \bfy\ or \bfY\ cannot
be adjacent to a layer labeled \bfp\ or \bfP.
\end{pro}

Next we prove ({\it Cf.}  \cite[Corollary~6.2]{rs1:1996}):

\begin{pro}
\label{pro:comp-free-iff-unlabelled} Let $R$ be a region of the
graphic and $(s,t) \in R$. Let $R'$ denote the image of $R$ under
$(s,t) \to (t,s)$. Then the intersection of $F_1(s,t)$ and
$F_2(s,t)$ is compression free and contains an essential curve if
and only if $R$ and $R'$ are both unlabeled.

Similarly, let $l$  be a layer and $t \in l$. Then
the intersection of $\s_t$ with $\gst$ is compression free yet
contains an essential curve if and only if $l$ is unlabeled.
\end{pro}

\begin{proof}
Suppose $R$ and $R'$ are both unlabeled.  By
Lemma~\ref{lem:lowercase} the absence of lowercase label in $R$
implies that $F_1(s,t) \cap F_{2}(s,t)$ contains a curve that is essential
in $F_{1}(s,t)$ or $F_{2}(s,t)$ (or both).  Let $\gamma$ be such a curve.
Assume (for contradiction) that $\gamma$ is essential in
$F_1(s,t)$ but not in $F_2(s,t)$ and denote by $D \subset
F_2(s,t)$ the disk $\gamma$ bounds in $F_2(s,t)$. Consider all
curves of $D \cap \s_t$ that are essential in $F_{1}(s,t)$. An innermost
such curve shows that $R$ has an uppercase label, contradiction.
Next assume (for contradiction) that $\gamma$ is essential
in $F_2(s,t)$ but not in $F_1(s,t)$. Then $f(\gamma)$ is
essential in $F_1(t,s)$ but not in $F_2(t,s)$ implying an upper
case label for $R'$, again contradicting our assumption. Hence the
absence of uppercase labels in $R$ and $R'$ implies that every
curve of $F_1(s,t) \cap F_2(s,t)$ is essential or inessential in both
surfaces, \ie, the intersection is compression free.

The converse is similar and we outline it here.
Suppose $R$ or $R'$ is labeled.  If $R$ has a lowercase
label then for $(s,t) \in R$, $F_1(s,t) \cap F_2(s,t)$ contains only
inessential curves by 
Lemma~\ref{lem:lowercase}.   If $R'$ has a lowercase label, then by the same
lemma, for $(s,t) \in R'$, $F_1(s,t) \cap F_2(s,t)$ contains only inessential
curves; since $F_1(t,s) = f(F_2(s,t))$ and $F_2(t,s) = f(F_1(s,t))$, applying
the involution we see that for $(s,t) \in R$,  $F_1(s,t) \cap F_2(s,t)$
contains only inessential curves as well. If $R$ has an uppercase label, then
by definition~\ref{dfn:labels} for $(s,t) \in R$ 
some curve of  $F_1(s,t) \cap F_2(s,t)$ is essential in
$F_1(s,t)$ and inessential in $F_2(s,t)$.  Finally, if $R'$ has an uppercase
label than for $(s,t) \in R'$,  some curve of  $F_1(s,t) \cap F_2(s,t)$ is
essential in $F_1(s,t)$ and inessential in $F_2(s,t)$.  Applying the
involution we see that for $(s,t) \in R$ some curve of  $F_1(s,t) \cap
F_2(s,t)$ is essential in $F_2(s,t)$ and inessential in $F_1(s,t)$.

The claim for layers follows from the claim for regions since every layer $l$
is contained in a region $R$ for which $R = R'$ and has the same labels as $R$.
\end{proof}

Let $l_p = (a',a)$ (for some $a,a' \in [-\infty,\infty]$) be the highest layer
labeled \bfp\ or \bfP. Since for $t << 0$ the label is \bfp\ the layer $l_p$
exists and since for $t >> 0$ the label is \bfy\ the layer $l_p$ is not the
topmost layer; hence $a \in \mathbb{R}$. Let $l_y = (b,b')$ (for some $b,b'
\in [-\infty,\infty]$)  be the first layer past $l_p$ labeled \bfy\ or \bfY.
Since the topmost layer is labeled \bfy\ the layer $l_y$ exists. By
Proposition~\ref{pro:adjacent-layers} the layers $l_p$ and $l_y$ cannot be
adjacent; hence $a<b$.  By choice of $l_p$, the layers between $l_p$ and $l_y$
are not labeled \bfp\ or \bfP, and by choice of $l_y$ they are not labeled
\bfy\ or \bfY.  Hence all the layers in $(a,b)$ are unlabeled and by
Proposition~\ref{pro:comp-free-iff-unlabelled} the corresponding surfaces have
compression free intersection, yet their intersection has an essential curve;
this completes the proof of Theorem~\ref{thm:compr-free}(1).

Let $t_0$ be a point $a < t_0 < b$ and suppose there is a region $R$ of the
Graphic adjacent to $(t_0,t_0)$ corresponding to an intersection which is
either 
not compression free, or consists entirely of inessential \scc s.  Since every
regular point $a < t < b$ is unlabeled, $(t_0,t_0)$ is not in the interior of
$R$; hence $(t_0,t_0)$  is a vertex of $R$. Let $R'$ be the image of $R$ under
$(s,t) \to (t,s)$ (note that $(t_0,t_0)$ is a vertex of $R'$ as well).   By
Proposition~\ref{pro:comp-free-iff-unlabelled}  either $R$ or $R'$ is labeled.
If the label at $R$ or $R'$ is \bfp\ or \bfP\ (resp. \bfy\ or \bfY) we shorten
the interval $(a,b)$ by replacing $a$ by $t_0$ (resp. replacing $b$ by
$t_0$). Repeating this process if necessary, we may assume all the regions
near every point of $(a,b)$ are unlabeled, some region near $a$ is labeled
\bfp\ or \bfP, and some region near $b$ is labeled \bfy\ or \bfY.

This completes the proof of Theorem~\ref{thm:compr-free}.
\end{proof}

Note that the proof gave us a little more than we bargained for: we have
control over the labels appearing near $a$ and $b$.  We need to improve the
intersection from ``compression free'' to ``essential''.  This is achieved in
the next section; for the remainder of this section we follow a technique of
\cite{rs1:1996} to  control inessential curves appearing in $(a,b)$.  Fix
$\epsilon>0$ small enough so that distance between any two critical points is
greater than $2\epsilon$.   Denote the critical points in $(a - \epsilon, b +
\epsilon)$ by $a=t_0 < t_1 < \cdots < t_{n-1} < t_n=b$.  Consider the
following embedding of the interval $(a-\epsilon,b+\epsilon)$ in the parameter
square, denoted $\delta$: send $t \in [a,b]$ to $(t,t)$, send $(a-\epsilon,a]$
into a region labeled \bfp\ or \bfP\ that is adjacent to $(a,a)$, and send
$[b,b+\epsilon)$ into a region labeled \bfy\ or \bfY.  More specifically, let
$R_1$ be the region containing $(a,t_1)$.  If $a$ is a single critical point
we embed $(a-\epsilon,a]$ in the diagonal.  If $a$ is a double  critical point
and there is a region labeled \bfp\ or \bfP\ that is adjacent to $R_1$ we
embed $(a-\epsilon,a]$ in that region.  Note that in that case moving from
$(a-\epsilon,a)$ to $(a,t_1)$  we cross only one of the critical points at
$a$ while tangent to the other, which we may ignore.  Finally, if $a$ is a
double critical point and both regions adjacent to $R_1$ are unlabeled  we
embed $(a-\epsilon,a]$ in the diagonal.  $(b,b+\epsilon)$ is treated
similarly.

$\delta$ gives an isotopy of two surfaces denoted $F_1(t)$ and $F_2(t)$
($(a-\epsilon < t < b+\epsilon)$).  We label points of $(a - \epsilon , b +
\epsilon)$ according to the intersection of $F_1(t)$ and $F_2(t)$, as in
Definition~\ref{dfn:labels}.

\begin{pro}
\label{pro:compr-free}
With the hypothesis of Theorem \ref{thm:compr-free} there are  families of
surfaces $F_1(t)$, $F_2(t)$ (for $t \in (a - \epsilon , b + \epsilon)$) with
the following properties: (1) $F_1(t)$ is isotopic to $\s_t$.\ (2) For $t \in
[a,b]$, $f(F_1(t)) = F_2(t)$.\  (3) $(a,b)$ has a neighborhood $N$ in the
parameter square so that every regular point in $N$ is unlabeled.  \ (4)
Arbitrarily close to $a$ (resp. $b$) there is a region labeled \bfp\ or \bfP\
(resp. \bfy\ or \bfY).

Furthermore, we may assume that for any regular point $t \in (a - \epsilon , b
+ \epsilon)$ the  intersection of $F_1(t)$ with $F_2(t)$ contains at most one
invariant inessential \scc\ or a pair of involute inessential \scc s.  When it
does, the layer containing $t$ is bounded by a single or double saddle on one
side and a single or double center on the other side, and the intersection in
the adjacent layers contains no inessential curves.
\end{pro}

\begin{proof}
The proof is based on \cite{rs1:1996}; we need to verify that it
works in the invariant setting.  We induct on the number of critical points in
$(a-\epsilon,b+\epsilon)$ that involve an inessential \scc.  Below, we modify
$F_1(t)$ and $F_2(t)$ by removing inessential \scc s via disk swaps or
introducing inessential \scc s via fixed or involute centers.  By definition
an inessential \scc\ bounds a disk on both $F_1$ and $F_2$  and therefore its
image is an inessential \scc\ as well; thus we may perform the disk swap
invariantly.  It is easy to see that this does not change labels.  We note
that we never change compressions or essential \scc s of $F_1 \cap F_2$, hence
there is a natural bijection between these curves before and after the
modification and we may talk of ``the same curves'' and ``the same saddles''.

We begin with $(a - \epsilon,a+\epsilon)$.  Assume first that the label at $a
- \epsilon$ is \bfp.   By isotopy (equivariant if $(a-\epsilon,a)$ was mapped
to the diagonal) we remove all inessential curves of $F_1(t) \cap F_2(t)$ for
$t \in (a-\epsilon,a+\epsilon)$.  At $a$ a single invariant inessential \scc\
(resp. two involute inessential \scc s) get pinched to form two or three
essential curves (resp. two essential curves).  We introduce a new critical
point at $a - \frac{1}{2}\epsilon$ as follows: if at $a$ two involute
inessential curves get pinched, the critical point at $a -
\frac{1}{2}\epsilon$ corresponds to involute  centers where the necessary pair
of curves appear.  If at $a$ a single invariant inessential curve gets
pinched we create this curve using an invariant center at $a -
\frac{1}{2}\epsilon$.  Note that this can be done: the invariant inessential
curve bounds two disks from $F_1$ and $F_2$ that (using an innermost disk
argument) we may assume are disjoint.  These disks bound an invariant ball that
by the Brouwer Fixed Point Theorem contains a fixed point
of $f$.  Isotoping $F_1(t)$ and $F_2(t)$ to that fixed
point we create an invariant center (we will often use this construction in
the proof of this proposition).  This defines the isotopy for $t \in (a -
\epsilon,a+\epsilon)$.  If the label at $(b,b+\epsilon)$ is \bfy, we modify
$(b - \epsilon,b + \epsilon)$ similarly.  Note that $(a - \epsilon,a +
\epsilon)$ and $(b - \epsilon,b + \epsilon)$ fulfill the requirements of
Proposition~\ref{pro:compr-free}. 

Next assume that  the label at $a-\epsilon$ is \bfP.  In that case, after
removing all inessential \scc s from the intersection at $(a-\epsilon,a)$ the
curves giving rise to  compressions still exist.  The critical point at $a$
is a (single or double) saddle that destroys these curves.  In the case of a
single saddle, this saddle cannot join an inessential \scc\ to a compression
or the labels would not change; therefore the (one or two) curves involved in
$a$ were not removed and we may cross $a$.  In the case of a double saddle
$(a-\epsilon,a+\epsilon)$ is embedded in the diagonal.   The two saddles
involve one, two, or three distinct curves.  If only one curve is involved it
is the compression which was not removed and we may cross $a$.  If two curves
are involved, at least one is a compression  and the other is the image of the
compression (note that a compression is never invariant) and hence is also a
compression; again we may cross $a$.  If three curves are involved, one  is
the compression.  Crossing only one of the two saddles we arrive at a region
adjacent to $R_a$, which by construction of the embedding of
$(a-\epsilon,b+\epsilon)$ we know is unlabeled.  Hence the curve attached to
the compression is not inessential, and similarly the third curve is not
inessential.  In this case too we may cross $a$.  In all cases, we constructed
the family $F_1(t), \ F_2(t)$ for $t \in (a-\epsilon, a +
\frac{1}{2}\epsilon)$.   If no inessential curves appeared at $a$ we extend
this family to $(a-\epsilon, a + \epsilon)$ without a change; otherwise, there
is either a single invariant inessential curve or a pair of involute
inessential curves.  We remove them using centers at $a+ \frac{1}{2}\epsilon$,
symmetrically to the construction of inessential curves described in the
previous paragraph.  This describes the modification at $(a -  \epsilon, a +
\epsilon)$ in this case.  If the label at $(b,b+\epsilon)$ is \bfY, we modify
$(b - \epsilon,b + \epsilon)$ similarly.  Note that $(a - \epsilon,a +
\epsilon)$ and $(b - \epsilon,b + \epsilon)$ fulfill the requirements of
Proposition~\ref{pro:compr-free}.  This concludes the base case of the
induction.

By the inductive hypotheses, we suppose the isotopy in
$(a-\epsilon,t_i+\epsilon)$ (for some critical $t_i \leq b$) fulfills the
requirements of Proposition~\ref{pro:compr-free}; moreover, by construction
$F_1(t_i+\epsilon) \cap F_2(t_i+\epsilon)$ contains no inessential curves.
The possibilities when passing from $t_0 - \epsilon$ to $t_0 +
\epsilon$ are (note that $t_i \in (a,b]$ hence every curve at $a < t < t_i$ is
either essential or inessential but not a compression): 

\begin{enumerate}
\item A single center (resp. double center) in which a one (resp. two)
  inessential \scc\ is created or destroyed. 
\item A single saddle (resp. double saddle) in which one (resp. two)
  inessential \scc\ is 
  attached to split off from another curve (which may or may not be
  inessential). 
\item A saddle or double saddle in which one inessential \scc\ becomes
  essential curves or compressions.
\item A double saddle in which two inessential \scc s become two
  essential curves or compressions.
\item All curves involved in $t_0$ are essential.
\end{enumerate}

In (1)--(2) the critical points are unnecessary, as they do not change
the pattern of essential curves.  Therefore we may ignore these critical
points and continue past $a$.  As above, in (3) we create the invariant \scc\
at $t_0 - \frac{1}{2} \epsilon$ and in (4) we create the two involute \scc s
at $t_0 - \frac{1}{2} \epsilon$.  In (5) there is nothing to do.

The surfaces $F_1(t_{i-1} + \epsilon)$, $F_2(t_{i-1} + \epsilon)$ and
$F_1(t_{i} - \epsilon)$, $F_2(t_{i} - \epsilon)$   are obtained from the
original surfaces by removing all inessential \scc s of intersection.  Hence
these surfaces are isotopic and we may extend the isotopy across $[t_{i-1} +
\epsilon,t_i - \epsilon]$.  Continuing in this  way we finally arrive at an
isotopy of $(a-\epsilon, t_{n-1} + \epsilon)$ that can be extended to
$(b-\epsilon,b+\epsilon)$ across $[t_{n-1} + \epsilon , b - \epsilon]$,
proving the proposition. 
\end{proof}

\begin{figure}
      \centerline{  \includegraphics[width=2in]{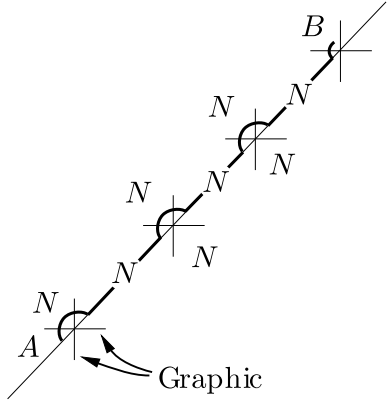}}
  \caption{The Generic Interval}
  \label{fig:isotopy}
  \end{figure}

For $(s,t) \in (a - \epsilon,b + \epsilon) \times (a - \epsilon,b + \epsilon)$
we construct the parameter square by setting $F_1(s,t) = F_1(s)$ and $F_2(s,t)
= F_2(t)$; the involution exchanges the surfaces along the diagonal only for
$t \in [a,b]$.  We perturb the parameter square fixing the diagonal to be
generic as we did in Section~\ref{sec:hhf}.   By construction $\delta$ is the
diagonal.  We perturb $\delta$ to obtain the {\it generic interval}: we move
$\delta$ slightly off $a$ and $b$ to be transverse to the Graphic; similarly,
near a double critical point (say $t_i$) we replace $(t_i - \epsilon, t_i +
\epsilon)$ (for some tiny $\epsilon$) by a small semicircle in $[a,b] \times
[a,b]$ that avoids the double point (see Figure \ref{fig:isotopy}).
By construction the generic interval is parameterized by $t \in
(a-\epsilon,b+\epsilon)$, starting at a region labeled \bfp\ or \bfP, going
through unlabeled regions to a region labeled \bfy\ or \bfY.

\section{Essential, Spinal Intersection}
\label{s:spinal}

\begin{dfns}
  \label{dfn:spinal-surf}
  \begin{enumerate}
    \item Let $S$ be a surface, and $K \subset S$ be an embedded graph.
    We say that $K$ {\it contains a spine} of $S$ if no component
    of $S$ cut open along $K$ contains a \scc\ that is essential
    in $S$.
    \item Let $F_1$, $F_2 \subset M$ be embedded surfaces, and
    $\Delta_2$ a set of compressing disks for $F_2$.  Suppose
    $F_1\cap F_2$, $F_1\cap \Delta_2$, and $F_1\cap \del \Delta_2$
    are all transverse.    We say that
    $F_1$ intersects $F_2 \cup \Delta_2$ {\it spinally}
    if $F_1 \cap (F_2 \cup \Delta_2)$ contains a spine of $F_1$.
    \item Let $F_1$, $F_2 \subset M$ be embedded surfaces.  We say
    that the intersection of $F_1$ and $F_2$ is {\it spinal} if
    there exists some set of compressing disks for $F_2$
    fulfilling condition (2) above or disks for $F_1$ fulfilling the same
    condition with the indices exchanged; for convenience we always
    assume that compressing disks are for $F_2$.\footnote{We often use (3)
    without any mention of the specific set of compressing disks.  We shall
    see that this makes sense, particularly in light of Lemma
    \ref{lem:m-cut-open} which shows the global nature of essential, spinal
    intersection.}
  \end{enumerate}
\end{dfns}

In this section we prove Theorem~\ref{thm:spinal}, which is a
combination of two theorems (one for a manifold admitting two
strongly irreducible Heegaard splittings and the other for a
manifold admitting a strongly irreducible Heegaard splitting and
an involution). 

Recall (Remark~\ref{rmk:non-haken}) that if $M$ is non-Haken then any minimal
genus \hhs\ is strongly irreducible; hence the theorem below is not vacuous:

\begin{thm}
\label{thm:spinal} 
Let $M$ be \assu\ manifold \hgt. Suppose that either $M$ admits
two strongly irreducible Heegaard surfaces $\s_1$ and $\s_2$ or a strongly
irreducible Heegaard surface $\s$ and an orientation preserving involution
$f$. Then we have:

\begin{enumerate}
  \item  $\s_1$ and $\s_2$ can be isotoped to intersect essentially and
  spinally.
  \item  $\s$ can be isotoped so that $\s$ and $f(\s)$ intersect essentially
  and spinally.
\end{enumerate}
\end{thm}

\begin{rmkk}
\label{rmk:improving-rs}{\rm
In \cite{rs1:1996} Rubinstein and Scharlemann prove a result very close to (1)
above: they show that $\s_1$ and $\s_2$ can be isotoped so that their
intersection is compression free, spinal, and contains at most one inessential
\scc.  (If we remove the inessential curve of intersection we may lose
spinality, so Theorem~\ref{thm:spinal}(1) does not follow.)  However, their
result is not quite strong enough for our purpose: in the next section we
prove that if $\s_1$ intersects $\s_2$ essentially and spinally then $M$ cut
open along $\s_1 \cup \s_2$ consists of handlebodies.  Existence of an
inessential curve of intersection allows for ``knotted handles'' and hence $M$
cut open along   $\s_1 \cup \s_2$ may not consists of handlebodies; it is
quite easy to  construct such examples.  }\end{rmkk}

\begin{proof}[Proof of Theorem \ref{thm:spinal}]
The method for finding a point that corresponds to spinal intersection is
given in \cite[Proposition 6.5]{rs1:1996}  where it is
shown that given an interval transverse to the Graphic, starting in a region
labeled \bfp\ or \bfP\ and ending in a region labeled \bfy\ or \bfY\ (such as
the generic interval constructed in the previous section)
there exists a set of compressing disks for one of the two surfaces (say
$\Delta_2$ for $F_2$) so that no component of $M$ cut open  along $F_2 \cup
\Delta_2$ is adjacent to itself,\footnote{We need this property for quoting
claims from \cite{rs1:1996} but we will not refer to it directly.} and for
some regular point $t$ in that interval the intersection of
$F_1(t)$ with $F_2(t) \cup \Delta_2(t)$ 
contains a spine of $F_1(t)$.  In this section, we show that
this point can be found on the diagonal and that the surfaces
corresponding to this point may be assumed to intersect essentially.  Note
that since the generic interval gives an ambient isotopy of $F_2$, it provides
an isotopy for $\Delta_2 = \Delta_2(t)$ as well ($(t)$ is suppressed
throughout this section).

We note that proving (1) of the theorem requires finding a point that
corresponds to spinal intersection in a layer that corresponds to essential
intersection, while (2) requires in addition that this point is on the
diagonal.  We will concentrate on (2) in this section and (1) will follow from
the argument here and the isotopy constructed in~\cite{rs1:1996} that has all
the properties of the generic interval.  From here on, we will not refer to
(1) directly.  

\begin{dfn}
\label{dfn:regula-2} Let $F_1$, $F_2$ and $\Delta_2$ be as above.
A point on the generic interval is called {\it regular} if the
intersections $F_1 \cap F_2$, $F_1 \cap \del \Delta_2$, and $F_1
\cap \Delta_2$ are all transverse, {\it critical} otherwise.  

After a
small perturbation of $\Delta_2$ (if necessary) we may assume there are only
finitely many critical points. The intervals obtained by cutting
the generic interval open along the critical points are called
{\it sublayers}.
\end{dfn}

The following lemma provides conditions to preserve spinality near
saddles.  A saddle move is similar to a boundary compression, and
crossing a saddle is equivalent to isotoping one of the surfaces
across a disk (say $\delta$) so that $\del \delta = (\delta \cap
F_1) \cup (\delta \cap F_2)$, where $\delta \cap F_1$ and $\delta
\cap F_2$ are two arcs meeting at their endpoints. We say that
$\delta$ {\it defines} the saddle. While the interior of $\delta$
is disjoint from $F_1$ and $F_2,$ it may intersect $\Delta_2$.  An
arc of $\Delta_2 \cap \delta$ has two endpoints, either both on
$F_1$, or one on $F_1$ and one on $F_2$, or both on $F_2$.  We say
that these arcs are of {\it type 1-1, 1-2, 2-2} (respectively).
See Figure \ref{fig:arc-types}.
\begin{figure}
\centerline{  \includegraphics[height=1.5in,width=3in]{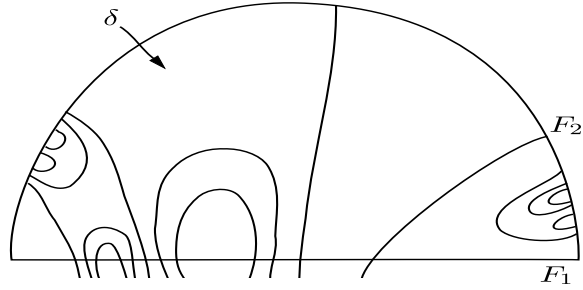}}
\caption{Arcs of $\Delta_2 \cap \delta$ (on $\delta$)}
\label{fig:arc-types}
\end{figure}

\begin{lem}
\label{lem:spinal-preserved-saddle}

Let $R_1$ and $R_2$ be
adjacent regions in the graphic, and suppose the critical point
separating $R_1$ from $R_2$ is a saddle.  Suppose further that
$\Delta_2 \cap \delta$ contains no type 1-1 arcs.

Then the intersection in $R_1$ is spinal if and only if the
intersection in $R_2$ is.
\end{lem}

\begin{proof}[Proof of Lemma \ref{lem:spinal-preserved-saddle}]
Assume the intersection is spinal in, say, $R_1$; we will show it is spinal in
$R_2$ as well.  Let $F_1$, $F_2$, and $\Delta_2$ be as above, intersecting
essentially.

First we show that after isotoping $\Delta_2$ if necessary we may
assume that $\Delta_2 \cap \delta$ contains no inessential \scc s:
let $\gamma \subset \Delta_2 \cap \delta $ be a \scc, chosen to be
innermost in $\delta$.  By isotopy of $\Delta_2$, we can replace
the disk $\gamma$ bounds in $\Delta_2$ by the disk it bounds in
$\delta$, and by small perturbation push this disk off $\delta$.
Since $\gamma$ was chosen innermost in $\delta$, $\Delta_2$
remains embedded. It is easy to see that $|\Delta_2 \cap \delta|$
was reduced by at least one; we need to show that the intersection
is still spinal.  The only change to $F_1 \cap (F_2 \cup
\Delta_2)$ is removing \scc s of $F_1 \cap \mbox{int}\Delta_2$.
None of these curves connects to any other component of $F_1 \cap
(F_2 \cup \Delta_2$). (We call such components {\it isolated \scc
s}.) If an isolated \scc\ were essential in $F_1$, then a parallel
copy of it would contradict spinality; hence isolated \scc s are
inessential in $F_1$ and removing them from $F_1 \cap (F_2 \cup
\Delta_2)$ does not change spinality.

Suppose that $\delta \cap \Delta_2 \neq \emptyset$.  Consider a  component
(say $T$) of $\delta$ cut open along $\Delta_2$ that contains a point of
$(\delta \cap F_1) \cap (\delta \cap F_2)$.  First assume $T$ contains no 2-2
arcs; recall that by assumption there are no 1-1 arcs.  Therefore $T$ is a
triangle with a single 1-2 arc on its boundary (see the leftmost triangle in
Figure~\ref{fig:arc-types}).
\begin{figure}
\centerline{  \includegraphics[width=4in]{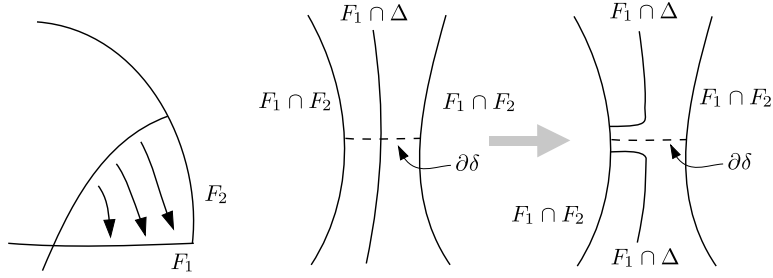}}
\caption{Removing arcs of type 1-2}
\label{fig:remiving-1-2-arcs}
\end{figure}
We use $T$ to guide an isotopy of $\Delta_2$ that removes the given 1-2 arc
from $\delta \cap \Delta_2$, see Figure~\ref{fig:remiving-1-2-arcs}.  Let
$\alpha$ be a \scc\ on $F_1$ disjoint from $F_2 \cup \Delta_2$ after the
isotopy.  It is easy to see that $\alpha$ is homotopic on $F_1$ to a curve
disjoint from $F_2 \cup \Delta_2$ before the isotopy
(Figure~\ref{fig:remiving-1-2-arcs} shows $F_1 \cap (F_2 \cup \Delta_2)$
``before and after''; $\alpha$ is not shown).  Since we assumed the
intersection to be spinal before the isotopy, $\alpha$ must be inessential on
$F_1$.  Hence the intersection is spinal after the isotopy.  This reduces
$|\delta \cap \Delta_2|$.

Next suppose that a component of $\delta$ cut open along
$\Delta_2$ that contains a point of $(\delta \cap F_1) \cap
(\delta \cap F_2)$ does contain 2-2 arcs;  we use $\delta$ to
guide an isotopy of $\Delta_2$ sliding these arcs off $\delta$. As
a result of this isotopy, for every 2-2 arc removed a pair of arcs
are added to $F_1 \cap (F_2 \cup \Delta_2)$ but nothing is
removed, hence the intersection is still spinal. This too reduces
$|\delta \cap \Delta_2|$.

Since $|\delta \cap \Delta_2|$ is being reduced this process must terminate;
when it does, $\delta \cap \Delta_2 = \emptyset$. We now cross the saddle.
After crossing the saddle the pattern of intersection between $F_1$ and $F_2
\cup \Delta_2$ changes only near the saddle point where two parallel arcs (say
horizontal) are replaced by two vertical arcs, denoted $v_1$ and $v_2$.
Clearly if we add an arc connecting $v_1$ to $v_2$ to $F_1 \cap (F_2
\cup \Delta_2)$ the intersection will become spinal.  We obtain this by the
move shown in Figure~\ref{fig:gecko-finger}.
\begin{figure}
      \centerline{  \includegraphics[width=4in]{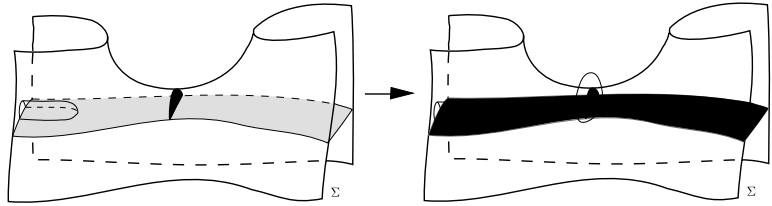}}
\caption{Making the intersection spinal after the saddle}
\label{fig:gecko-finger}
\end{figure}
It is described as follows: after crossing the saddle we obtain a disk similar
to $\delta$ defining the same saddle from the opposite side; denote this
disk $\delta^{-1}$.  Then $\delta^{-1} \cap F_1$ is an arc connecting  $v_1$
to $v_2$.  Let $V$ be the component of $M$ cut open along $F_2$ 
containing $\delta^{-1}$. Denote the frontier of a neighborhood 
of $\delta^{-1}$ in $V$ by
$D$.\footnote{In other words, $D$ is the disk obtained from two parallel
copies of $\delta^{-1}$ connected together on the other side of $F_1$ so that
$\del D \subset F_2$.    Note that $D$ is a boundary parallel disk  in $V$.}
$D \cap F_1$ consists of exactly two arcs parallel to $\delta^{-1} \cap
F_1$. Let 
$\alpha$ be an arc on $F_2$ connecting $\del \Delta_2$ to $\del D$, missing
$\del \Delta_2$ in its interior (\ie, an outermost arc; $\alpha$ may intersect
$F_1$ in its interior).  Band-connect sum a disk of $\Delta_2$ to $D$ along
$\alpha$; this 
changes $\Delta_2$ by an isotopy.  After the band sum, $F_1 \cap (F_2 \cup
\Delta_2)$ consists of the intersection prior to the band sum union an arc for
every point of $\alpha \cap F_1$ union $D \cap F_1$. It is now clear that $F_1
\cap (F_2 \cup \Delta_2)$ contains a spine of $F_1$, proving
Lemma~\ref{lem:spinal-preserved-saddle}.
\end{proof}

It follows immediately from Definition~\ref{dfn:labels} (labels) that at a
region labeled \bfp\ or \bfP\ there is a meridian disk for $F_1(t)$ that is
purple near its boundary (\ie, ``below'' $F_1(t)$) and intersects $F_1(t)$ in
curves that are all  inessential in $F_1(t)$ ({\it cf.}
\cite[Section~8]{rs1:1996}). Similarly, in a layer labeled \bfy\ or \bfY\
there exists a meridian disk for $F_1(t)$ that is yellow near its boundary
(\ie, ``above'' $F_1(t)$).  This motivates the following labeling scheme for
sublayers: 

\begin{dfn}
\label{dfn:second-labels}
The label \bel\ is used in a sublayer where there is a compressing disk for
$F_1$ that is disjoint from $F_2 \cup \Delta_2$, is below $F_1$, and 
intersects $F_1(t)$ (if at all) only in curves that are inessential in
$F_1(t)$.  The label \ab\ is defined similarly.
\end{dfn}

Strong irreducibility implies that each sublayer has at most one label and
adjacent sublayers cannot be labeled \bel\ and \ab. Since the generic interval
starts with a layer labeled \bfp\ or \bfP\ (and hence with a sublayer labeled
\bel) and ends with a layer labeled \bfy\ or \bfY\ (and hence with a sublayer
labeled \ab) some sublayer is unlabeled. From \cite[Section~8]{rs1:1996} we
have:

\begin{pro}[\cite{rs1:1996}]
\label{pro:spinal-iff-unlabeled} 
Let $t$ be a regular point.  $F_1(t) \cap (F_2(t)  \cup \Delta_2(t))$ contains
a spine of $F_1(t)$  if and only if $t$ is in an unlabeled sublayer.
\end{pro}

Thus, the generic interval described in Proposition~\ref{pro:compr-free}
contains a point $t$ that corresponds to spinal intersection. However, $t$ may
not have all the properties required by Theorem~\ref{thm:spinal}, specifically:

\begin{enumerate}
\item $t$ may be off the diagonal and is separated from it by centers.
\item $t$ may be off the diagonal and is separated from it by saddles.
\item $t$ may be in a sublayer where one or two inessential \scc s of
  intersection exist.
\end{enumerate}

Everything we said until this point is true for {\it any} isotopy of
$\Delta_2(t)$.  We exploit this flexibility and design an isotopy of
$\Delta_2(t)$ that helps us deal with the three problems listed above, then
combine the three cases to prove the theorem.

\noindent Case 1: The unlabeled sublayer is separated from
the diagonal by centers, denoted $b$ and $c$ (say $b < c$).  We may assume
that crossing $b$ from left to right an inessential curve appears (otherwise,
we reverse $t$).  Then by construction crossing $c$ from left to right another
inessential  curve appears.  Let $a < b$ and $d > c$ be points of the generic
interval on the diagonal and $\epsilon > 0$ small enough so that the only
critical points of $F_1(t) \cap F_2(t)$ in $(a - \epsilon , d + \epsilon)$ are
$b$ and $c$.

Fix $t_0 \in (b,c)$.  Then there is an arc $\alpha_c$ so that one endpoint of
$\alpha_c$ is $\alpha_c \cap F_1(t_0)$, the other is $\alpha_c \cap F_2(t_0)$,
and crossing $c$ is equivalent to isotoping $F_1(t_0)$ along $\alpha_c$ and
pushing a small disk of $F_1(t_0)$ across $F_2(t_0)$.  We say that $\alpha_c$
{\it defines the center} $c$.  We change the isotopy of $\Delta_2(t)$ as
follows:  we reparameterize $\Delta_2(t)$ in $(a-\epsilon,b)$ so that all the
critical points of $F_1(t) \cap \Delta_2$ and $F_1(t) \cap \del \Delta_2$ are
in $(a - \epsilon , a)$.  In $[b,t_0]$ $F_1 \cap \Delta_2$  and $F_1 \cap \del
\Delta_2$ have no critical points.  In $(t_0,c)$ we slide $\Delta_2$ off
$\alpha_c$.  Thus $F_1 \cap \Delta_2$ has exactly $|\alpha_c \cap \Delta_2|$
critical points, each introducing an isolated \scc\ to $F_1(t) \cap (F_2(t)
\cup  \Delta_2(t))$.  $F_1 \cap \del \Delta_2$ has no critical point in $(t_0
, c)$.    In $[c , d + \frac{1}{2} \epsilon]$ there are no critical points of
$F_1(t) \cap \Delta_2(t)$ or $F_1(t) \cap \del \Delta_2$.  In $(d +
\frac{1}{2} \epsilon, d + \epsilon)$ we isotope $\Delta_2$ to its original
configuration.    For a regular value $t \in (a,d)$ the difference between
$F_1(t) \cap (F_2(t) \cup \Delta_2(t))$ and $F_1(a) \cap (F_2(a) \cup
\Delta_2(a))$ is isolated curves; hence the intersection is spinal in $t$ if
and only if it is spinal in $a$, and we may assume Case~(1) does not happen.
Moreover, we have control over the labels of sublayer: either all sublayers of
$(a,d)$ are unlabeled (if the intersection is spinal) or all are labeled, and
since  adjacent labeled sublayers have the same label we conclude that {\it
either both $a$ and $d$  are both unlabeled or both are labeled and the labels
at $a$ and $d$ are the same.}

\medskip

\noindent Case 2: The unlabeled sublayer is separated from the
diagonal by saddles.  Similar to Case~(1) denote the saddles $b < c$ and let
$a<b$ and $d>c$ be points of the generic interval on the diagonal, and
$\epsilon > 0$ small enough so that the only critical point of $F_1(t) \cap
F_2(t)$ in $(a - \epsilon , d + \epsilon)$ are $b$ and $c$.   At $a$ there
exist two disks $\delta_b, \ \delta_c$ defining the saddle  $b,\ c$,
respectively.  (Recall the construction of $\delta$ in the paragraph
preceeding Lemma~\ref{lem:spinal-preserved-saddle}.)  Since moving along the
diagonal both saddles are crossed simultaneously $\del \delta_b \cap \del
\delta_c = \emptyset$, and applying a standard innermost disk argument we may
assume that $\delta_b \cap \delta_c = \emptyset$.
Similar to the proof of Lemma~\ref{lem:spinal-preserved-saddle} we use
$\delta_b$ and $\delta_c$ to guide an isotopy of $\Delta_2$ off $\delta_b$ and
$\delta_c$. (This changes the sublayers, and since $\Delta_1 \cap \delta_b$
and $\Delta_2 \cap \delta_c$ may have 1-1 arcs we cannot assume the labels do
not change.)  Isotope $\Delta_1$ in $(a,d)$ so that $F_1 \cap \Delta_2$ and
$F_1 \cap \del \Delta_2$ have no critical points in $[a,d]$.  (Thus $(a,b)$
and $(c,d)$ are contained in one sublayer each, and $(b,c)$ is a sublayer).
In $(d , d + \epsilon)$ isotope $\Delta_2$ to its original configuration.  By
Lemma~\ref{lem:spinal-preserved-saddle} if the sublayer $(b,c)$ is unlabeled
so is the sublayer containing $(c,d)$.  Hence, if a regular value $t \in
(b,c)$ corresponds to spinal intersection so does the regular value $d$ on the
diagonal.  As in Case~(1) we have a little more:
Lemma~\ref{lem:spinal-preserved-saddle} implies that the sublayer containing
$(a,b)$ is unlabeled if and only if the sublayer $(b,c)$ is.  As adjacent
sublayers have the same labels we again conclude that {\it either both $a$ and 
$d$  are unlabeled or both are labeled and the labels at $a$ and $d$ are 
the  same.}

\medskip

\noindent Case 3.  The unlabeled sublayer corresponds to an intersection that
contains one or two inessential \scc s.   Let $l$ be a layer containing
inessential curves.  By Proposition~\ref{pro:compr-free} $l$ is bounded on one
side by a (single or double) center and the other side by a (single or double)
saddle.  Say the center is at $c$ and the saddle at $s$.  For convenience we
assume $c < s$ (the other case is similar), so $l = (c,s)$.  Let $\epsilon >
0$ be small enough so that $c$ and $s$ are the only critical points on the
generic interval in $(c - \epsilon, s + \epsilon)$.   For convenience we
assume the semicircles of the generic interval have radius
$\frac{1}{3}\epsilon$.  For $t \in (c,s)$ there are one or two $\delta$ disks
that define the saddles; if  there are two $\delta$ disks we may assume (as in
Case~(2)) that they are disjoint.  In $(c-\epsilon, c)$ there are one or two
$\alpha$ arcs that define the centers.  By the construction in
Proposition~\ref{pro:compr-free} the saddles at $s$ involve the inessential
curves of $(c,s)$; hence the $\delta$ disks cannot be seen in $(c - \epsilon,
c)$.  However, for $t \in (c - \epsilon, c - \frac{1}{3}\epsilon)$ we can find
the trace  of the  $\delta$ disks as  disks are disjointly embedded in their
interior but not on their boundary.  The boundary of each disk  consists of
four arcs, one on $F_1(t)$, one on $F_2(t)$, and between them two arcs on the
$\alpha$ arcs.  There are three cases, but their treatment is identical: in
the case of a single center and a single saddle we see a single $\delta$ disk
attached to itself along the single $\alpha$ arc to form an annulus
$\mathcal{A}$.  In case of a single center and two saddles we see two $\delta$
disks attached to each other along a single $\alpha$ arc, each disk forming an
embedded annulus (say $A_1, \ A_2$), with $A_1 \cap A_2 = \alpha$.  We take
$\mathcal{A}$ to be $A_1 \cup A_2$.  The case of a double center and a single
saddle is impossible since there would still be inessential curves in $(s , s
+ \epsilon)$, contradicting Proposition~\ref{pro:compr-free}.  Finally, in
case of double center and double saddle, the two $\delta$ disks are glued to
each other along the two $\alpha$ arcs forming a single annulus $\mathcal{A}$.
In all three cases we see an annulus-like complex $\mathcal{A}$ which is
homeomorphic to either a circle cross an interval or a wedge of two circles
cross an interval, and the arcs $\alpha$ are contained in $\mathcal{A}$ and
have the form one or two points cross interval.
 
We are now ready to describe the isotopy of $\Delta_2$:  in $(c - \epsilon, c
-\frac{2}{3}\epsilon)$ slide $\Delta_2$ off the $\alpha$ arcs.  As before this
creates isolated \scc s.  At $c - \frac{2}{3} \epsilon$ the arcs of $\Delta_2
\cap  \mathcal{A}$ come in three flavors, arcs of types $1-1$, $1-2$, and
$2-2$, where an arc is of type $i-j$ if it has one boundary component on $F_i$
and the other on $F_j$.  In $(c - \frac{2}{3} \epsilon,c - \frac{1}{3}
\epsilon)$  we use $\mathcal{A}$ to guide an isotopy of $\Delta_2$ that
removes all $1-1$ arcs (the so-called karate-chop).  After crossing the
centers, the $\delta$ disks contain no $1-1$ arcs.  We isotope $\Delta_2$ in
$(c + \frac{1}{3}\epsilon , s - \frac{1}{3}\epsilon)$ to remove the $1-2$ and
$2-2$ arcs.  After crossing the saddles near $s$, we isotope $\Delta_2$ in $(s
+ \frac{1}{3}\epsilon , s + \epsilon)$ to its original configuration.  After
this isotopy, if some layer in the semicircle $(c - \frac{1}{3}, c +
\frac{1}{3} \epsilon)$ is unlabeled then so is the layer past $c + \frac{1}{3}
\epsilon$, as addition of isolated curves at the centers cannot change
spinality.  If some label in $(c + \frac{1}{3} \epsilon,  s + \frac{2}{3}
\epsilon)$ is unlabeled then by Proposition~\ref{lem:spinal-preserved-saddle}
the region containing $s + \frac{1}{3} \epsilon$ is unlabeled.  We conclude
that {\it if some layer in $(s - \epsilon, c + \epsilon)$ is unlabeled than
some layer in $(c - \epsilon,c - \frac{1}{3} \epsilon)$ or in $(s +\frac{1}{3}
, s + \epsilon)$ is unlabeled.}

For proving Theorem~\ref{thm:spinal}(1): in \cite{rs1:1996} Rubinstein and
Scharlemann give an isotopy of $\s_1$ and $\s_2$ with the properties listed in
Proposition~\ref{pro:compr-free} (with no reference to invariance, of course).
Theorem~\ref{thm:spinal}(1) follows from that and the argument in Case~(3)
above.

We combine the three cases to prove Theorem~\ref{thm:spinal}(2): starting with
the generic interval $(a-\epsilon,b+\epsilon)$, we isotope $\Delta_2$ in a
neighborhood of any layer that contains an inessential curve as described in
Case~(3) above.  Next, given a double critical point not on the boundary of an
layer containing inessential curves (say $t_0$), we isotope $\Delta_2$ near it
as described in Cases~(1) or (2) above.  The generic interval starts at a
sublayer labeled \bel\ and ends at a sublayer labeled \ab\ and is transverse
to the Graphic; by \cite{rs1:1996} some sublayer of the generic interval is
unlabeled, and by Cases~(1), (2), and (3) above there exists a point on the
diagonal corresponding to essential, spinal intersection.

This completes the proof of Theorem~\ref{thm:spinal}.
\end{proof}

\section{$M$ cut open along $\s \cup f(\s)$}
\label{sec:m-cut-open}

This section is devoted to the proof of Theorem~\ref{thm:m-cut-open}. 
The proofs of cases~(1) and (2) are identical.  For simplicity we use the
notation $\s_1$ and $\s_2$ in the proof, (2) follows by setting $\s = \s_1$
and $f(\s) = \s_2$.    In Theorem~\ref{thm:spinal} we established the
existence of an isotopy of  $\s_1$ and $\s_2$ so that the intersection of
$\s_1$ and $\s_2$ is essential and spinal.  Theorem~\ref{thm:m-cut-open}
follows from that and the following lemma that originally appeared
in \cite{rieck-proc}.  For completeness we bring it here with its
proof.

\begin{lem}

\label{lem:m-cut-open} Let $\s_1$ and $\s_2$ be \hhs s
intersecting spinally and essentially.  Then the components of $M$
cut open along $\s_1 \cup \s_2$ are handlebodies.
\end{lem}

\begin{proof}

Since the intersection is spinal there exists a complete set of
compressing disks $\Delta_2$ for one of the surfaces (say $\s_2$) so
that $\s_1 \cap (\s_2 \cup \Delta_2)$ contains a spine of $\s_1$. By
definition of spinal intersection, $\s_1$ is incompressible
in the complement of $\s_2 \cup \Delta_2$.  (Note that components of $\s_1$
cut open along $\s_2 \cup \Delta_2$ may compress, but any curve of
$\s_1$ cut open along $\s_2 \cup \Delta_2$ that is compressed is
inessential in $\s_1$.)

We may assume that $\s_1 \cap \Delta_2$ consists of arcs only: let
$\gamma$ be a \scc\ in $\s_1 \cap \Delta_2$.  Since the intersection
is spinal, $\gamma$ bounds a disk in $\s_1$.  Passing to an
innermost such, we see a disk whose interior intersects neither
$\Delta_2$ nor $\s_2$ (by essentiality). We now use this disk to
isotope $\Delta_2$ and reduce $|\Delta_2 \cap \s_1|$.

Let $B$ be some component of $M$ cut open along $\s_2 \cup
\Delta_2$, and $c$ some component of $\s_1 \cap B$.  We show that
$c$ is a disk.  Assume for contradiction $c$ is not a disk.  Since the
intersection is spinal, every curve on $c$ is trivial in $\s_1$.  Hence $c$ is
a punctured disk. Let
$\gamma$ be one of the punctures, and $D\subset \s_1$ the disk it
bounds (see Figure~\ref{fig:c-not-disk}).
\begin{figure}
\centerline{  \includegraphics[height=1.5in,width=2.5in]{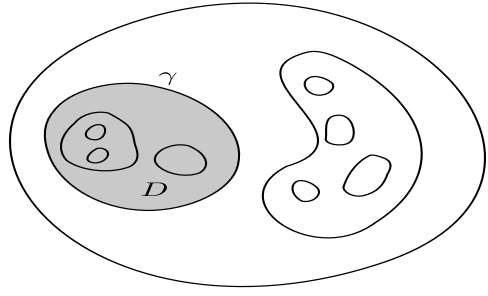}}
\caption{When $c$ is not a disk} \label{fig:c-not-disk}
\end{figure}
By assumption $\del D = \gamma \subset \del B,$ and $N_D(\del(D)) \cap
B=\gamma$ (that is, near its boundary $D$ is outside $B$).
Since the intersection of $\s_1$ and
$\s_2$ is essential $\gamma \not\subset \s_2$.
Since $\s_1 \cap \Delta_2$ consists of arcs, $\gamma \not\subset
\Delta_2$.
Hence $\gamma$ must have parts on $\s_2$ and parts on $\Delta_2$ (say
$\Delta_2$ is {\it above} $\s_2$).  Clearly part of $D$ is below $\s_2$. But
the boundary of this part of $D$ is a non-empty collection of \scc s
in $\s_1\cap \s_2$, all inessential in $\s_1$, contradicting
essentiality.

$M$ cut open along $\s_2 \cup \Delta_2$ consists of balls.
Since the pieces of $\s_1$ in each of these balls
are disks, they further chop these balls up
into balls, that is to say, $M$ cut open along $\s_1 \cup \s_2 \cup \Delta_2$
consists of balls.  As we saw, $\s_1 \cap \Delta_2$ consists entirely of arcs
and therefore $\Delta_2$ cut open along $\s_1$ consists of disks.
Attaching the balls of $M$ cut open along $\s_1 \cup \s_2 \cup \Delta_2$ to
each other via these disks we get handlebodies of $M$ cut open along $\s_1
\cup \s_2$.
\end{proof}

\section{The an-annular complex $\ca$}
\label{sec:modify-ca}

Using $\s$ found in Theorem~\ref{thm:m-cut-open} we define $\ca$ to be $\s
\cup f(\s)$. $\ca$ is a complex, mostly a surface, but with some points that
are not surface points.  At these points $\ca$ looks like the intersection of
two surfaces.  Denote this set by \sca.  However, in this section we will
modify \ca\ and it will no longer be the union of two surface; we think
of \ca\ as a collection of embedded surfaces with boundary, disjoint in their
interiors, and with images of any two boundary components either disjoint or
equal.  Then \sca\ is the union of boundary components.  Denote the genus of
$\s$ by $g$.  \ca\ has the following properties:

\begin{props}
  \label{prop-ca:prime}
\begin{description}
\item[A] $\chi(\mbox{\rm\ca})\geq 4-4g$.
\item[B] All components of $M$ cut open along {\rm\ca} are handlebodies.
\item[C'] No piece of $\mbox{\rm\ca\ }\setminus\mbox{\rm\sca}$ is a disk.
\item[D'] Every curve of {\rm \sca} is the union of an even number of boundary
components.
\item[E] Every torus embedded in {\rm\ca} bounds a solid torus.
\item[F] {\rm\ca} is invariant under the involution.
\end{description}
\end{props}

Properties {\bf  A, D'} and {\bf F} are obvious.  Properties {\bf B} and {\bf
C'} are Theorem \ref{thm:m-cut-open}. Property {\bf E} was proved by Kobayashi
and Rieck in \cite[Corollary~1.3]{rieck-kobayashi.1}.

However, Property {\bf C'} is insufficient as components of $\mbox{\rm\ca\
}\setminus\mbox{\rm\sca}$ may be annuli, preventing an Euler characteristic
count. We need to replace Properties {\bf C'} and {\bf D'} with a stronger
version, Properties {\bf C} and {\bf D} below.  Achieving these properties is
the context of this section and requires us to modify \ca. To see the relation
between Properties {\bf C'} and {\bf D'} and Properties {\bf C} and {\bf D} we
mention that in the process of modifying \ca, we remove from \ca\ a
neighborhood of \sca\ and replace it by the boundary of that neighborhood, so
all curves of \sca\ are arranged along tori and have valence three, where the
{\it valence} of a curve of \sca\ is the number of surfaces adjacent to it
locally.  (This does not completely describes the modification we perform.)

The closure of a component of $\ca \setminus \sca$ is called a {\it sheet}.
In Property~{\bf A} stated below we also consider a manifold admitting two
strongly irreducible Heegaard splittings of genera $g_1$ and $g_2$.

\begin{props}
  \label{prop-ca}
  \begin{description}
  \item[A] $\chi(\mbox{\rm\ca})\geq 4-4g$ (or $4-2(g_1+g_2)$).
  \item[B] All components of $M$ cut open along {\rm\ca} are handlebodies.
  \item[C] Every curve of {\rm \sca} is the union of three boundary
    components, one of a sheet with negative Euler characteristic and two of
    annular sheets.  These annuli close up, together with other annular
    sheets, to form tori bounding solid tori (denoted $\{V_i\}_{i=1}^n$).  For
    each $i$, $V_i \cap \mbox{{\rm\ca}} = \emptyset$.
  \item[D] For each $i$, the number of annuli forming $\del V_i$ is even.
  \item[E] Every torus embedded in {\rm\ca} bounds a solid torus.
  \item[F] {\rm\ca} is invariant under the involution.
  \end{description}
\end{props}

\begin{example} {\rm
\label{ex:no-annuli-in-RS} Whenever $\s \cup f(\s)$ contains no
annuli, removing from \ca\ a neighborhood of \sca\ and replacing it by the
boundary of that neighborhood is sufficient for achieving
Properties~\ref{prop-ca}.  The following example shows that this requirement
is sometimes impossible to impose: let $M$ be a genus 2 manifold admitting a
free involution and let $\s$ be a genus two \hhs\ for $M$.  Suppose $\s$
intersects $f(\s)$ essentially and spinally and without annuli.  It is easy to
see that $\s$ cut open along $\s \cap f(\s)$ consists of two components,
either both once punctured tori or both pairs of pants (similarly, $f(\s)$ cut
open along $\s \cap f(\s)$ consists two components homeomorphic to the
components $\s$ cut open along $\s \cap f(\s)$).  Denote by $V$ one of the
handlebodies obtained by cutting $M$ open along $\s$.  We then see that $\del
(V \cap f(V))$ is a genus two surface and using Lemma~\ref{lem:m-cut-open} we
deduce that $V \cap f(V)$ is a genus two handlebody. Therefore $f|_{V \cap
f(V)}$ is a free involution, and the quotient of $V \cap f(V)$ by the
involution has Euler characteristic $-\frac{1}{2}$, contradiction.  }
\end{example}

We now state the main theorem of this section.  In this theorem, $\s$ (resp.
$\s_1$ and $\s_2$) are the surfaces found in Theorem~\ref{thm:spinal} for a
manifold with involution (resp. a manifold admitting two strongly irreducible
Heegaard splittings).

Recall (Remark~\ref{rmk:non-haken}) that if $M$ is non-Haken then any minimal
genus \hhs\ is strongly irreducible; hence the theorem below is not vacuous:

\begin{thm}
\label{thm:complex}
Let $M$ be \assu\  manifold \hgt\ admitting an \op\ involution $f:M \to M$ and
a strongly irreducible \hhs\ $\s$ of genus $g$ (resp. two strongly irreducible
Heegaard surfaces or genera $g_1$ and $g_2$).

Then there exists a complex $\mbox{\rm\ca}\subset M$ fulfilling Properties
\ref{prop-ca} (resp. Properties~\ref{prop-ca}{\bf A}--{\bf E}).  Moreover, if
$\s$ (resp. $\s_1$ and $\s_2$) is the surface found in
Theorem~\ref{thm:spinal}, we may assume that $(\s \cup f(\s)) \setminus
(\cup_{i=1}^n V_i) = \ca \cap (\cup_{i=1}^n V_i)$ (resp. $(\s_1 \cup \s_2)
\setminus (\cup_{i=1}^n V_i) = \ca \cap (\cup_{i=1}^n V_i)$).
\end{thm}

\begin{rmkk}
\label{rmk:no-change-outside-msts}{\rm

The proof is constructive, giving an algorithm that takes the surface $\s$
(resp. $\s_1$ and $\s_2$) found in Theorem~\ref{thm:m-cut-open} as input and
starting with $\ca = \s \cup f(\s)$ (resp. $\s_1 \cup \s_2$) modifies \ca\ in
finitely many steps until arriving at a complex (still denoted \ca) fulfilling
Properties~\ref{prop-ca}. (The algorithm given here is an equivariant version
of the algorithm given in \cite{rieck-proc}.)  }
\end{rmkk}

\begin{question}
\label{que:link-exterior}{\rm
We view the cores of the solid tori described in Property~\ref{prop-ca}{\bf C}
as a link in $M$.  Not every link can arise this way (for example, one can see
that Properties~\ref{prop-ca} imply an upper bound on the Heegaard genus of
the link exterior and the number of its components).  We ask what (and how
many) links arise in this way and what other properties do they have.}
\end{question}

\begin{proof}[Proof of Theorem~\ref{thm:complex}]
The proofs for manifold with involution and manifolds containing two strongly
irreducible Heegaard surfaces are identical except for the invariance
requirement (Property~{\bf F}), which makes the latter strictly easier.  We
therefore concentrate on the former only.

Starting with $\ca\ = \s \cup f(\s)$, we modify \ca\ to fulfill Properties
\ref{prop-ca}. As noted above Properties~\ref{prop-ca:prime} are already
satisfied, and (unless replaced with stronger properties)  they must be
preserved throughout the work; that is to say they are {\it invariants} of the
algorithm.

During the modifications of \ca\ we construct the solid tori $V_i$ and enlarge
them, step by step.  However, we never modify \ca\ outside these solid tori.
This guarantees that for the final complex $(\s \cup f(\s)) \cap (M \setminus
(\cup_{i=1}^n V_n)) = \ca \cap (M \setminus (\cup_{i=1}^n V_n))$ as required.
We will not refer to this again.

The next invariant counts the number of sheets attached to a solid torus $V
\subset M$ with $\del V \subset \ca$.  We count with multiplicity, that is, a
sheet with $n$ boundary components on a $V$ is counted as $n$ sheets attached
to $V$.  (Sheets inside $V$ are not counted.)  The invariant is:

\begin{inv}
\label{inv:even-number-attached-to-V}
Let $V \subset M$ be a solid torus, $\del V \subset\ca$.  Then the number of
sheets attached to $V$ is even.
\end{inv}

\begin{proof}[Proof of Invariant~\ref{inv:even-number-attached-to-V}]

Let $\gamma \subset \del V$ be a curve of \sca.  Since the valence of $\gamma$ is
four and exactly two sheets attached to $\gamma$ are part of $\del V$, at
$\gamma$ there are either zero, one, or two sheets attached to $V$. We need to
show that the number of the curves with one sheet attached to $V$ is even.  By
Property~{\bf C'} $\gamma$ is essential in $\del V$.

Assume first that the slope defined by \sca\ is the meridian of $V$.  We can
then remove $\del V$ from \ca, obtaining an immersed surface $\mathcal{S}$.
If at $\gamma$ two sheets of $\mathcal{S}$ are attached to $\del V$ from
outside (resp. inside) $V$, push $\mathcal{S}$ near $\gamma$ out of
(resp. into) $V$, removing $\gamma$ from $\del V$.  If some component (say
$F$) of $\mathcal{S} \cap int(V)$ is not a meridian disk then $F$ is either
boundary parallel, compressible, or boundary compressible.  In the first case
$F$ is a boundary parallel annulus (since $\del F$ is essential in $\del V$)
and we push $F$ out of $V$ without changing the parity of the number of
sheets attached to $V$.  In the second case, we compress $F$.  If $F$ is
boundary compressible but not boundary parallel then $F$ is compressible, so
we may ignore the third case (see, for example,
\cite[Lemma~2.7]{rieck-kobayashi.1}).  Finally we see that every component of
$\mathcal{S}$ in $\mbox{int}(V)$ is a meridian disk; by construction the
number of meridian disks of $\mathcal{S} \cap V$ has the same parity as the
number of curves on $\del V$ where a single sheet was attached to $V$.  We
constructed $\mathcal{S}$ by removing $\del V$ from \ca, isotopy and
compression.  Hence $\mathcal{S}$ is homologous to the null-homologous complex
\ca\ and the number of times $S$ intersects the core of $V$ is even.  This
number is exactly the number of meridian disks of $\mathcal{S} \cap V$,
proving Invariant~\ref{inv:even-number-attached-to-V} in this case.

Next assume that the slope defined by \sca\ is not meridional.  If at $\gamma$
zero (resp. one, two) sheets are attached to $\del V$ from outside $V$, then
two (resp. one, zero) sheets are attached to $\del V$ from {\it inside} $V$.
Hence the number of sheets attached to $V$ (from outside) is even if and only
if the number of sheets attached to $\del V$ from inside is even.  Let $F$ be
the closure of a component of $\s \cap \mbox{int}(V)$ ($F$ may intersect
$f(\s)$ in its interior and so may not be a sheet).  In \cite[Section 2]{mr}
and \cite[Theorem~3.3]{MR99h:57040} it was shown that if a strongly
irreducible Heegaard surface $\s$ intersects a solid torus $V$ so that each
curve of $\s \cap \del V$ is a non-meridional essential curve of $\del V,$
then a component $F$ of $\s \cap V$ is either an annulus, or a twice punctured
torus, or a four times punctured sphere; in particular $|\del F|=2$ or $|\del
F|=4$.  The same holds for every component of $f(\s) \cap V$.  Summing up
these numbers gives the number of sheets attached to $\del V$ from inside;
hence this number is even as required.
\end{proof}

The main tool used in this section is:

\begin{dfn}
\label{dfn:mst} A solid torus $V$ embedded in $M$ is called a {\it
\mst} if $\del V \subset \ca$ and $V$ is maximal with respect to inclusion
among all such solid tori.
\end{dfn}

By definition a \mst\ is an embedded solid torus; in particular, a solid torus
embedded in its interior but not in its boundary cannot be a \mst.  Let
$\{V_i\}_{i=1}^n$ be the set of all \msts\ in $M$, which is finite since the
complex \ca\ is.

We would like \msts\ to be disjoint; this is not quite the case.  For future
reference we state this lemma for any complex \ca\ fulfilling Property~{\bf
E}; in particular, Property~{\bf F} (invariance) is not used in the proof.

\begin{lem}
\label{lem:msts-intersection} Let {\rm\ca} be a complex fulfilling Property~{\bf
  E}.  Then any two distinct \msts\ are either disjoint or intersect in a
single \scc\ that is essential in the boundary of both and longitudinal in (at
least) one.
\end{lem}

\begin{proof}[Proof of Lemma~\ref{lem:msts-intersection}]

Let $V_1$ and $V_2$ be distinct \msts\ so that $V_1 \cap V_2 \neq
\emptyset$. We first show that $V_1 \cap V_2$ is a single \scc.  Let $W = V_1
\cup V_2$. Let $\{N_i\}_{i=1}^k$ be the closures of the components of $M
\setminus W$.  If, for some $i$, $\del N_i$ contains an embedded surface (say
$S$) then $S$ is a torus (it has zero Euler characteristic since it is made up
of annuli, and is orientable since it locally separates $W$ from $N_i$ in the
orientable manifold $M$). By Property {\bf E}, $S$ bounds a solid torus in
$M$, and by maximality this solid torus cannot contain $V_1$ or
$V_2$. Therefore it must contain $N_i$ and we conclude that (since $N_i$ is
connected) the solid torus is $N_i$ itself.  If $N_i$ is a \st\ for all $i$
then $M$ is the union of $N(W)$ with solid tori.  This gives a decomposition
of $M$ into solid tori that intersect in annuli.  If some slope is meridional,
$M$ is reducible or a lens space; else, $M$ is a \sfs; all conclusions
contradict our assumptions.

Therefore we may assume that some component (say $N_1$) is not a solid torus
and hence no component of $\del N_1$ is an embedded surface.  Thus there
is some curve on  $\del N_1$ (say $\gamma$) where $V_1$ is tangent to $V_2$. A
neighborhood of $\gamma$ in \ca\ separates a neighborhood of $\gamma$ in $M$
into four regions, two non-adjacent (say east and west) from $V_1$ and $V_2$,
and the other two (north and south) from $N_1$. If $V_1 \cap V_2 = \gamma$ we
are done.  Thus we may assume $\del V_1 \cap \del V_2$ contains at least one
more component.  Note that $\del V_1$ is a torus, formed by gluing an annulus
connecting (say) the southeast corner of $\gamma$ to the northeast corner to
itself along $\gamma$.   Denote this annulus by $A_{V_1}$, and similarly
denote $A_{V_2}$ the annulus connecting the southwest corner to the northwest
corner, so that gluing $A_{V_2}$ to itself at $\gamma$ gives $\del V_2$.  By
assumption $A_{V_1}$ is not disjoint from $A_{V_2}$ in its interior.  Let $T'$
be an embedded torus obtained from cut and pasting annuli of $A_{V_1}$ and
$A_{V_2}$ cut open along $A_{V_1} \cap A_{V_2}$.  Then $T'$ is a toral
component of $\del N_1$ and by the previous paragraph $N_1$ is a torus,
contradicting out assumption.  This shows that $V_1 \cap V_2$ is a single
curve.

Next we show that the slope of  $V_1 \cap V_2$ is longitudinal in $V_1$ or
$V_2$.  For contradiction assume that the slope of the intersection is not
longitudinal in either \st.  If it is meridional in both then $M$ contains a
non-separating sphere and if it is meridional in one and cabled in the other
(\ie, neither meridional nor longitudinal) then $M$ contains a lens space
summand, both contradicting our assumptions.   So we may assume the slope is
cabled in both. Consider $W$ be $N(V_1 \cup V_2)$ which is a \sfs\ over $D^2$
with exactly two exceptional fibers.  Denote $\del W$ by $T$.  If $T$ bounds a
\st\ then either $M$ reduces or $M$ is a \sfs.  Thus $T$ is a torus not
bounding a solid torus.  By assumption $M$ is irreducible and a-toroidal and
therefore $T$ bounds a knot exterior contained in a ball, say $X$ (for details
see, for example, \cite{rieck-kobayashi.1}).  If $X$ were $\mbox{cl}(M
\setminus W)$ then $T$ would be essential, contradicting our
assumptions. Hence $X = W$.  In \cite[Theorem~1.1]{rieck-kobayashi.1}
Kobayashi and Rieck proved that if a strongly irreducible Heegaard surface
intersects a torus bounding a knot exterior contained in a ball in curves that
are all essential in the torus, then the slope defined by these curves is
meridional.  In our case $T \cap \s$ is the slope of a regular fiber in the
Seifert fibration by construction, which is not meridional (note that  $X$ is
a torus knot exterior), contradiction.
\end{proof}

\noindent We now modify \ca\ in four steps (we do not rename \ca\ after each
step):

\medskip

\noindent {\bf Step One: amalgamating \msts.}

\begin{dfn}
\label{dfn:amalgamation} Let $V_1, \ V_2 \subset M$ be solid tori
such that $\del V_1, \ \del V_2 \subset \ca$ and $V_1 \cap V_2$ is a \scc\
$\gamma$, so that $\gamma$ is essential in $\del V_1$ and $\del V_2$ and
longitudinal in at least one of $V_1$, $V_2$.  Let $N(\gamma)$ be a small
neighborhood of $\gamma$, invariant if $\gamma$ is.  Replacing \ca\ by $(\ca
\setminus \ca \cap N(\gamma)) \cup (\mbox{cl}(\del N(\gamma) \setminus (V_1
\cup V_2)))$ is called {\it amalgamating} $V_1$ and $V_2$ along $\gamma$ (or
simply amalgamating along $\gamma$, or amalgamating $V_1$ and $V_2$).  The two
annuli $\mbox{cl}(\del N(\gamma) \setminus (V_1 \cup V_2))$ are denoted $A_1$
and $A_2$, the solid torus obtained by amalgamating $V_1$ and $V_2$ is denoted
$V,$ and its boundary is denoted $T$. Note that $V_1, \ V_2 \subset V$ and
exactly one curve was removed from \sca; no other curve of \sca\ has changed.
\end{dfn}

Suppose there exist \msts\ (say $V_1$ and $V_2$) so that $V_1 \cap V_2 \neq
\emptyset$. Amalgamate $V_1$ and $V_2$ (which can be done by
Lemma~\ref{lem:msts-intersection}).  We show that the resulting solid torus
$V$ is a \mst: let $U$ be a \mst\ containing $V$.  If $U$ is embedded prior to
the amalgamation then $V_1, \ V_2 \subset U$, contradicting their maximality.
Else, prior to the amalgamation $U$ is pinched at $\gamma$ and broken up to
two solid tori, one containing $V_1$ and the other containing $V_2$.  By
maximality, these solid tori are $V_1$ and $V_2$ themselves and $U = V$.
Therefore $V$ is a \mst\ as desired.  We verify Property~{\bf E}:

\begin{lem}
\label{lem:prop-E-after-amalgamation}
Let $V_1$ and $V_2$ be \msts\ in a complex fulfilling Property~{\bf E} and
assume $V_1$ can be amalgamated with $V_2$.  Then {\rm\ca} fulfills Property~{\bf
E} after amalgamation.
\end{lem}

\begin{proof}
Let $T\subset \ca$ be a torus after the amalgamation.  Then one of the
following holds:

\begin{enumerate}
\item $A_1 \not\subset T$ or $A_2 \not \subset T$.
\item $A_1 \subset T$ and $A_2 \subset T$.
\end{enumerate}

In Case~(1) $T$ is embedded in \ca\ prior to amalgamation.  Since \ca\
fulfills Property~{\bf E} before the amalgamation $T$ bounds a solid torus. In
Case~(2), prior to the amalgamation there are two tori (say $T'$ and $T''$) so
that $T' \cap T'' = \gamma$ and $T$  is obtained from $T'$ and  $T''$ via
surgery.  By property~{\bf E}, $T'$ and $T''$ bound solid tori (say $V'$ and
$V''$ respectively; note that if $V' \subset V''$ we cannot amalgamate the
two).  Let $U'$, $U''$ be the \msts\ containing $V'$, $V''$ respectively. Then
$\gamma \subset U'$ and hence so are at least two of the four sheets adjacent
to $\gamma$.  Thus $U' \cap V_1$ or $U' \cap V_2$ contains a sheet, and by
Lemma~\ref{lem:msts-intersection} either $U' = V_1$ or $U' = V_2$, say the
former.  Similarly either $U'' = V_1$ or $U'' = V_2$.  Since $T' \cap T'' =
\gamma$ we  see that $U'' = V_2$. Therefore $V' \cap V'' \subset U' \cap U'' =
V_1 \cap V_2 = \gamma$, and $V'$  can be amalgamated to $V''$ along $\gamma$.
Clearly, $T$ bounds the  amalgamation of $V'$ and $V''$.
\end{proof}

If $\gamma$ is an invariant curve, we perform the amalgamation invariantly.
Else, we amalgamate along $f(\gamma)$; we verify that this can be done: Let
$V_3$ be a \mst\ distinct from $V_1, \ V_2$ above.  If $V$ (the result of
amalgamating $V_1$ and $V_2$) intersects $V_3$, by
Lemma~\ref{lem:msts-intersection} the intersection is a single essential curve
that is longitudinal in at least one of the two solid tori.  Thus we can
amalgamate along $f(\gamma)$ (either amalgamating $V$ and $V_3$ or
amalgamating two \msts, both distinct from $V$).  After this, Property~{\bf F}
is recovered.

We continue amalgamating as long as possible, always performing the
amalgamation invariantly.  This process reduces $|\sca|$ and hence terminates.
When it does, any two \msts\ are disjoint and \ca\ is invariant. We may now
replace Lemma~\ref{lem:msts-intersection} with the stronger property below,
which is our next invariant:

\begin{inv}
\label{inv:disjoint-msts}
Any two \msts\ are disjoint.
\end{inv}

We check invariants:

\begin{description}
\item[Property A]  $\chi(\ca)$ has not changed.

\item[Property B]  The new components of $M$ cut open along \ca\ are solid
  tori. 

\item[Property C']  The new sheets are annuli.

\item[Property D']  Some curves are removed from \sca\ and the number of
  sheets attached to all other curves is unchanged.

\item[Property E]  See Lemma~\ref{lem:prop-E-after-amalgamation}

\item[Property F] By construction.

\item[Invariant~\ref{inv:even-number-attached-to-V}] In the Proof of
Lemma~\ref{lem:prop-E-after-amalgamation} we saw that any new solid torus
(after amalgamating $V_1$ and $V_2$) is the amalgamation of two solid tori
$V', \ V''$ at $\gamma$.  Prior to the amalgamation, the number of sheets
attached to $V', \ V''$ is even, and the number of sheets attached to the
amalgamation of $V'$ and $V''$  is the sum of these numbers minus four.
\end{description}

\medskip

\noindent {\bf Step Two: cleaning \msts.}  We remove from \ca\ every sheet
that is in the interior of a \mst.\footnote{Note that $\s \cup \gst$ may
contains many components inside a \mst.  In that case \ca\ will be modified
very drastically in Step Two.  For example, if $V$ is a \mst\ and $\ca \cap V$
looks like a grid cross $S^1$ then in Step Two many annuli are removed, which
is the reason this step is important for the algorithm constructed here.}   As
a result, the valence of curves of \sca\ on the boundary of each \mst\ is
either three or four.  Let $\gamma$ be a curve on the boundary of a \mst\ $V$
with valence four.  We equivariantly deform \ca\ by adding a small
neighborhood of $\gamma$ to $V$, splitting $\gamma$ into two curves of valence
three.  This completely describes the modification of \ca\ in Step Two.

We show that the tori embedded in \ca\ after Step Two are exactly the
boundaries of \msts\ before Step Two.  In one direction, if $V$ is a \mst\
prior to Step Two then clearly $\del V$ is a torus embedded in \ca\ after Step
Two.  For the other direction, let $T \subset \ca$ be an embedded torus after
Step Two.  Let $T'$ be the image of the embedding prior to Step Two.  If $T'$
is not embedded then $T'$ has a double curve (say $\gamma'$) on a valence four
curves of \sca\ on the boundary of a \mst.  Locally near $\gamma'$, $T'$ has
four annuli, two on $\del V$ and two attached to $\del V$.  Let $A'$ be one of
the annuli of $T'$ attached to $\del V$.  It is easy to use the annuli of $T'$
cut open along double curves to cut and paste an embedded torus (say $T''$)
containing $A'$.  By Property~{\bf E}, $T''$ bounds a solid torus and this
solid torus is contained in a \mst, say $V''$.  By construction $A' \subset
V''$ and therefore $V'' \neq V$ and $V'' \cap V \neq \emptyset$, contradicting
Invariant~\ref{inv:disjoint-msts}.

So we may assume that $T'$ is embedded.  Then by Property~{\bf E}  $T'$ bounds
a \st, say $V$.  $V$ is contained in some \mst, say $U$.  If $V \neq U$ then
parts of $\del V$ are in the interior of $U$ and are thrown out in Step Two,
contradicting choice of $T$.  Hence $V = U$ and $T = \del U$ as required.

This proves the following invariant, which is stronger than Property~{\bf E}
and therefore replaces it:

\begin{inv}
\label{inv:step-two}
Every torus embedded in \ca\ bounds a maximal solid torus that does not
intersect \ca\ in its interior.
\end{inv}

We call a curve $\gamma \in \sca$ that is on the boundary of a \mst\ a {\it
  boundary curve} and a sheet on the boundary of a \mst\ a {\it boundary
  sheet}.  If $\gamma \in \sca$ is a boundary curve then by
  Invariant~\ref{inv:disjoint-msts} it is on the boundary of exactly one \mst\
  and hence of the three sheets attached to $\gamma$ exactly two are boundary
  sheets.  We replace Property~{\bf D'} by Property~{\bf D''} to accommodate
  boundary curves:

\begin{description}
\item[Property D''] If $\gamma \subset \sca$ is not a boundary curve then
$\gamma$ is the union of four boundary components. Every boundary curve has
valence three.  The boundary of a \mst\ consists of an even number of boundary
sheets.
\end{description}

Note that Property~{\bf D''} implies
Invariant~\ref{inv:even-number-attached-to-V} and hence replaces it.  We now
check our invariants, proving Property~{\bf D''}.

\begin{description}
\item[Property A]  Since no sheet is a disk the Euler characteristic is no
  more negative than it was.
\item[Property B]  The new components of $M$ cut open along \ca\ are solid
  tori.
\item[Property C']  The only new sheets are boundary sheets, and they are all
  annuli.
\item[Property D''] For non-boundary curves there is nothing new to prove.
  For $\gamma \in \ca$ a boundary curve, this follows immediately from
  Invariant~\ref{inv:even-number-attached-to-V}.
\item[Property F]  Since the image of a \mst\ is a \mst, \ca\ is invariant.
\item[Invariant~\ref{inv:disjoint-msts}]  The set of \msts\ was not changed in
  Step Two.
\end{description}

\medskip

\noindent {\bf Step Three: curves of sing(\ca) not on \msts.} Let $\gamma$
be a curve of sing(\ca) not on the boundary of a \mst. Note that such a curve
was not changed from the original complex.

\begin{dfn}
A map from a torus into \ca\ is called {\it admissible} if it is a
homeomorphism on the torus except at a finite set of double curves. On double
curves the map is 2-to-1 into curves of \sca.
\end{dfn}

Thus each double curve either double covers its image or is identified with
another double curve. Since $\ca \setminus \sca$ contains no disks, the 
image of the torus
cut open along the double curves consists of annuli.  Note that we do not
require annuli adjacent to double curves to cross each other, that is to say,
an annulus coming from the south may be connected to an annulus from the east,
while an annulus from the north is connected to an annulus from the west (so
an admissible map need not be self-transverse as a map into $M$).

\begin{lem}
\label{lem:admis-maps} The only admissible maps are boundaries of
\msts.
\end{lem}

\begin{proof}
For contradiction assume that there exists an admissible map $g:T \to \ca$
that is not the boundary of a \mst. By Invariant~\ref{inv:step-two} every
torus embedded in \ca\ is the boundary of a \mst; therefore the map considered
is not an embedding and has a double curve in \sca, say $\gamma$.  Since
boundary curves have valence three, $\gamma$ is not a boundary curve.
Therefore $\ca\cap N(\gamma)$ was not changed in Steps One and Two, and $\ca
\cap N(\gamma)$ is the intersection of two  annuli. Thus $\gamma$ is the image
of two distinct curves on $T$ and these curves cut $T$ into two annuli, say
$A$ and $A'$.  Since both boundary components of $A$ map to $\gamma$, $A$
defines an admissible map with $\gamma$ on its boundary and fewer double
curves than $g$. Continuing in this way, we construct an embedding of the
torus into \ca\ that intersects some curve of \sca\ that is not on the
boundary of a \mst, contradicting Invariant~\ref{inv:step-two}. Thus every
admissible map is an embedding and hence the boundary of a \mst.
\end{proof}

Let $\gamma$ be a curve of \sca\ not on the boundary of a \mst.  Replace \ca\
by $(\ca\ \setminus N(\gamma)) \cup (\del N(\gamma))$, introducing a new solid
torus.  This construction can be done equivariantly by either considering
pairs of involute curves, or using the Invariant Neighborhood Theorem on
invariant curves.

Let $T$ be a torus embedded in \ca\ after Step Three.  It is easy to see that
$T$ has one of the following two forms: either prior to Step Three there is
some 
non-boundary curve $\gamma$ and $T$ is $\del N(\gamma)$, or prior to Step
Three $T$ is an admissible map.  Hence by Lemma~\ref{lem:admis-maps}, either
$T$ bounds a solid torus $V$ given by $N(\gamma)$ (for some non-boundary curve
$\gamma$) or $T$ bounds a solid torus $V$ that was a \msts\ prior to Step
Three.  Thus we see that after Step Three the set of \msts\ consists of
neighborhoods of non-boundary curves and \msts\ prior to Step Three; clearly,
distinct \msts\ are disjoint.  Every curve of \sca\ is a boundary
curve and every non-boundary sheet has its boundary on \msts.  We emphasize
that since maximal solid tori are disjoint, every curve of \sca\ is the
boundary of exactly two boundary sheets and one non-boundary sheet (although
all three may be annuli). We replace Properties~{\bf C'} and {\bf D''} by
Properties~{\bf C''} and {\bf D}, which are very close to the required
Properties~{\bf C}  and {\bf D}.  (In fact, if we could replace
``non-positive" in Property~{\bf C''} by ``negative" we'd be done.)

\medskip

\noindent {\bf Property C'':} Every curve of {\rm \sca} is the union of three
boundary components, one of a non-boundary sheet with non-positive Euler
characteristic and two annular boundary sheets.  These boundary sheets close
up, together with other boundary sheets, to form tori bounding solid tori.
These solid tori do not intersect \ca\ in their interior.

\medskip

\noindent {\bf Property D:}   The number of annuli forming each torus
described in Property~{\bf C} is even.

\medskip

We now check invariance of the properties achieved so far.

\noindent {\bf Property A:} The Euler characteristic was not changed in Step
Three.

\noindent {\bf Property B:}  All new components of $M$ cut open along \ca\ are
solid tori.

\noindent {\bf Property C'':} By construction.

\noindent {\bf Property D:}  Every new torus has four annuli on its boundary.

\noindent {\bf Property F:} By construction.

\noindent {\bf Invariant~\ref{inv:disjoint-msts}:} By construction.

\noindent {\bf Invariant~\ref{inv:step-two}:} By Lemma~\ref{lem:admis-maps}
and the construction.

\begin{rmk}
We pause for a moment to review what we achieved so far.  Recall from
Example~\ref{ex:no-annuli-in-RS} that at the onset our only concern were
annular sheets (of course now these sheets are best described as annular
non-boundary sheets). Many annular sheets were removed in Step Two.   The
crucial property we achieved by using \msts\ is that chains of (boundary and
non-boundary) annular sheets do not close up to form tori, except for boundary
of \msts.  This allows us to remove annular non-boundary sheets in Step Four.
\end{rmk}

\medskip

\noindent {\bf Step Four: getting rid of annular non-boundary
sheets.} Let $A$ be an annular non-boundary sheet. Assume (for contradiction)
that $A$ connects a \mst\ (say $V$) to itself.  We use $A$ and an annulus of
$\del V$ cut open along $\del A$ to form a torus, say $T$.  Let $U$ be the
\mst\ that $T$ bounds (which exists since every torus bounds a \mst).
Clearly, $U$ and $V$ are distinct \msts\ and $U \cap V \neq \emptyset$,
contradiction. Thus $A$ connects two distinct \msts, say $V_1$ and $V_2$.
Assume (for contradiction) that the slopes defined by $\del A$ on $\del V_1$
and $\del V_2$ are both not longitudinal.  If both slopes are  meridional then
$M$ contains a non-separating sphere and if one slope is meridional and the
other cabled then $M$ contains a lens space summand, both conclusions
contradicting our assumptions.  Finally, if the slope is cabled in both $V_1$
and $V_2$ then $N(V_1 \cup A \cup V_2)$ is a \sfs\ over the disk with exactly
two 
exceptional fibers.  As in the proof of Lemma~\ref{lem:msts-intersection} it
is easy to argue that $\del N(V_1 \cup A \cup V_2)$ is a torus not bounding a
solid torus.  Hence $\del N(V_1 \cup A \cup V_2)$  bounds a knot exterior $X$,
and (since $M$ is a-toroidal) $X = N(V_1 \cup A \cup V_2)$.  Then we have:

\begin{clm}
$\s \cap \del X$ or $f(\s) \cap \del X$ is non-empty and consists of fibers in
the Seifert fibration of $X$.
\end{clm}

\begin{proof}

Since $M$ cut open along \ca\ consists of handlebodies (and not compression
bodies) \ca\ is connected.  Since $\mbox{cl}(M \setminus X)$ is not a solid
torus, \ca\ is not contained in $X$.  Hence $\ca \cap \del X \neq \emptyset$.
A fiber in the Seifert fibration is given by a curve on $\del X$ parallel to
$\del A$; it is now easy to see that all curves of $\ca \cap \del X$ are
parallel (in $\del X$) to such a curve, and hence are fibers.

Denote the set of \msts\ $V_1,\dots,V_n$; by
Remark~\ref{rmk:no-change-outside-msts} $ (\s \cup f(\s)) \cap (M \setminus
(\cup_{i=1}^n V_n)) = \ca \cap (M \setminus (\cup_{i=1}^n V_n))$.  Since $\del
X \subset (M \setminus (\cup_{i=1}^n V_n))$, we have that $\s \cap \del X$ or
$f(\s) \cap \del X$ is non-empty and consists of fibers, proving the claim.
\end{proof}

However, by \cite[Theorem~1.1]{rieck-kobayashi.1} if a strongly irreducible
Heegaard surface $\s$ (or $f(\s)$) intersects a knot exterior $X$ contained in a
ball and $\s \cap \del X$ consists of a non-empty collection of curves
that are all essential in $\del X$ then these curves are meridional. The
meridian of a torus knot exterior is not a fiber, contradiction. We conclude
that $A$ connects two distinct \msts\ and is longitudinal in at least one of
the two.  Similar to Step One we amalgamate $V_1$ with $V_2$ along $A$ by
replacing \ca\ with $(\ca \setminus N(A)) \cup \mbox{cl}(\del N(A) \setminus
(V_1 \cup V_2))$, see Figure \ref{fig:step-four}.  Denote the new component of
$M$ cut open along \ca\ by $V$.   $V$ is a solid torus.
\begin{figure}
  \centerline{\includegraphics[width=1.5in]{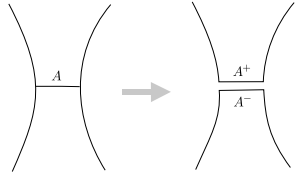}}
  \caption{Amalgamation along an annulus}
  \label{fig:step-four}
\end{figure}

Let $T$ be a torus embedded in \ca\ after the amalgamation.  Denote the two
parallel copies of $A$ in $\del N(A)$ by $A^+$ and $A^-$. First ({\it cf.}
Lemma~\ref{lem:prop-E-after-amalgamation}(1)) suppose  $A^+ \not\subset T$ or
$A^- \not\subset T$.  Then $T$ is embedded prior to the amalgamation and hence
$T$ is the boundary of a \mst\ (say $U$) prior to amalgamation.  It is
straight forward to see that $U \neq V_1$ and $U \neq V_2$, hence $U \cap V_1
= \emptyset$ and $U \cap V_2 = \emptyset$, that is, $U$ existed as a \mst\
prior to the amalgamation. Next, ({\it cf.}
Lemma~\ref{lem:prop-E-after-amalgamation}(2)) suppose $A^+ \subset T$ and $A^-
\subset T$. Denote  the boundary components of $A^+$ and $A^-$ by $S_1^+$,
$S_2^+$, $S_1^-$, and $S_2^-$ where $S_i^\epsilon = A^\epsilon \cap V_i\;
(\epsilon=\pm,\,i=1,2)$.  We follow $T$ along, starting at $S_1^+$, moving
across $A^+$ to $S_2^+$, continuing until we get to $A^-$.  Assume (for
contradiction) that the boundary component we get to is $S_1^-$.  Then the
annulus we have traversed from $S_1^+$ to $S_1^-$ forms an embedded torus (say
$T'$) prior to amalgamation.  Then $T'$ bounds a \mst\ (say $U'$) and it is
easy to see that $U' \neq V_1$, and $U' \cap V_1 \neq \emptyset$;
contradiction.  Thus we conclude that the boundary component we get to is
$S_2^-$.  The annulus we traversed from $S_2^+$ to $S_2^-$ forms an embedded
torus (say $T'$) prior to the amalgamation, and the \mst\ it bounds intersects
$V_2$.  Hence, this \mst\ is $V_2$. Similarly around $V_1$, and we conclude
that $T$ is simply the boundary of the new \st. Thus, the set of \msts\ after
amalgamation is exactly the set of \msts\ prior to amalgamation (except $V_1$
and $V_2$), with the new solid torus $V$ replacing $V_1$ and $V_2$.  Note that
distinct \msts\ are disjoint.

Exactly as in Step One we notice that any two \msts\ can be amalgamated: let
$A'$ be an annulus connecting \msts.  The argument in the beginning of Step
Four shows that $A$ connects distinct \msts\ and is longitudinal in at least
one.  If $A$ is not invariant we retrieve invariance of \ca\ by amalgamating
along $f(A)$.  Iterating Step Four as long as we can, the process reduces
$|\sca|$ and terminates when no non-boundary sheet is an annulus.  Note that
after Step Four, every torus embedded in \ca\ bounds a \mst\ and distinct
\msts\ are disjoint.

We now verify Properties~\ref{prop-ca}:

  \begin{description}
  \item[A] The Euler characteristic of \ca\ was not changed in Step Four.
  \item[B] All new components of $M$ cut open along {\rm\ca} are solid tori.
  \item[C]  By construction, every curve of \sca\ is on three sheets, two
  boundary sheets and one non-boundary sheet.  Boundary sheets close up to form
  tori bounding solid tori.  In Step Four we removed all non-boundary annular
  sheets.
  \item[D]  This property is preserved since after amalgamating (say) $V_1$
  and $V_2$ to obtain $V$, the number of sheets attached to $V$ is the sum of
  the sheets attached to $V_1$ and $V_2$ minus four.
  \item[E] By construction, every torus embedded in \ca\ bounds a (maximal)
  solid torus.
  \item[F] \ca\ is invariant by construction.
  \end{description}

This completes the proof of Theorem \ref{thm:complex}.
\end{proof}

\section{Constructing the invariant \hh\ surface $S$}
\label{sec:get-s}

\begin{proof}[Proof of Theorem~\ref{thm:first}]
By Theorem~\ref{thm:complex} $M$ admits a complex \ca\ fulfilling
Properties~\ref{prop-ca}{\bf A}--{\bf F}.  Recall that the solid tori
components of $M$ cut open along \ca\ are denoted $\{ V_i \}_{i=1}^n$.  Let
$\ka^*$ be the complex $\ca\ \cup (\cup_{i=1}^n V_i)$.  Note that $\chi(\ka^*)
= \chi(\ca)$ and therefore by Property~\ref{prop-ca}{\bf A} $\chi(\ka^*) \geq
4-4g$. By Property~{\bf F}, $\ka^*$ is invariant under $f$.  We call the
components of $\ka^* \setminus (\cup_{i=1}^n V_i)$ {\it sheets} (recall that
$V_i$ are close solid tori and therefore the annuli forming $\del V_i$ are not
sheets).

Let \ka\ be the complex obtained from $\ka^*$ by puncturing every sheet once
or twice (if necessary for invariance).  Note that the only sheets that are
punctured twice are sheets that are invariant but admit no fixed point, and
every such sheet has Euler characteristic divisible by 2.   By Property~{\bf
C} every sheet has negative  Euler characteristic, we see that every sheet
that is punctured once has Euler characteristic at least as negative as -1 and
every sheet that is punctured twice has Euler characteristic at least as
negative as -2.  Every puncture reduces the Euler characteristic by exactly
one, and we see that the Euler characteristic is doubled at worst, that is to
say, $\chi(\ka) \geq 8-8g(\s)$.

Let $S = \del N\ka$.  Since \ka\ is invariant, so is $S$. On one side (away
from \ka) $S$ bounds components  of $M$ cut open along \ca\ (that are all
handlebodies by Property~{\bf  B}) glued to each other along  disks that
correspond to the punctures of \ka.  Thus $S$ bounds a handlebody on that
side.  On the side containing \ka, $S$ bounds the solid tori $V_i$ glued along
pieces of the form punctured sheet cross interval.  Since a punctured sheet
deformation retracts to a spine that contains the boundary of the sheet, a
punctured sheet cross interval deformation retracts to a neighborhood of that
spine.  It is now easy to see that this component of $M$ cut open along $S$ is
obtained from $\cup_{i=1}^n V_i$ by attaching 1-handles; hence it too is a
handlebody.  Thus $S$ is an invariant Heegaard surface and the complementary
handlebodies are invariant.  We calculate  its genus: $\chi(S) =  \chi(\del
N\ka) = 2 \chi(N\ka) = 2\chi(\ka) \geq 16-16g(\s)$.  Thus $2-2g(S) \geq 16 -
16 g(\s)$.  Solving for $g(S)$ we get $g(S) \leq \bddup$.  This completes the
proof of Theorem \ref{thm:first}.
\end{proof}

\section{Constructing a \hhs\ for $N$}
\label{sec:hhs-for-quotient}

\begin{proof}[Proof of theorem \ref{thm:main}]
The surface $S$ found in Theorem \ref{thm:first} is an invariant \hhs\ for $M$
and the involution preserves the sides of $S$. Pick a handlebody of
$M\setminus S$, say $H$.

\begin{clm}
\label{clm:quotient}
The quotient of $H$ by the involution is a handlebody.
\end{clm}

We prove the claim by induction on the genus of $H$, denoted $g(H)$.  For
balls, this is a result of Waldhausen~\cite{MR0236916}.  Assume $g(H)>0$.  By
the Equivariant Disk Theorem $H$ admits equivariant essential disks, either an invariant
disk $D$ or two disjoint disks $D, \ D'$ that are involutes of each other.
The image of the equivariant disks is a single disk $f(D)$.  We cut $H$ along
the equivariant disks, obtaining $H_D$.  $H_D$ consists of one, two or three
handlebodies, all of genus lower than $g(H)$.  We cut $f(H)$ along $f(D)$
obtaining $f(H)_{f(D)}$.  The projection $f$ induces a cover $f|{H_D}: H_D \to
f(H)_{f(D)}$.  By induction, $f(H)_{f(D)}$ consists of handlebodies.  Gluing
these handlebodies to each other along $f(D)$ we see that the image of $H$ is
a handlebody.  This proves the claim.

We see that $N$ cut open along the image of $S$ (denoted $S/(f)$) consists of
two handlebodies, and therefore $S/(f)$ is a  Heegaard surface for $M/(f)$.
In Section~\ref{sec:get-s} we saw that $\chi(S) \geq 16-16g(\s)$.  If $f|_S$
has no fixed points then $f$ induces an unbranched cover $F|_S : S \to S/(f)$.
In that case the Euler characteristic is multiplicative and we get: $2 -
2g(S/(f)) = \chi(S/(f)) \geq 8 - 8g(\s)$; solving for $g(S/(f))$ we see that
$g(S/(f)) \leq \bdddown.$  It is easy to see that if the cover $f|_{S} : S \to
S/(f)$ is branched the genus of $S/(f)$ is even lower.

This completes the proof of Theorem \ref{thm:main}.
\end{proof}

\section{Bounding the Bridge Number of the Branch Set}
\label{sec:bridge}

In this section we prove Theorem~\ref{thm:bridge}, bounding the complexity of
the branch set of the double cover $f:M \to N$.  The branch set is a link in
$N$, denoted $L$.  To measure the complexity of $L \subset N$ we fix a
Heegaard surface $F$ for $N$ and isotope the link to intersect each of the
handlebodies of $N$ cut open along $F$ in boundary parallel arcs.  (To see
that this is possible, pick any Heegaard function corresponding to $F$ and
pull the maxima of $L$ above zero and the minima below.)  We define:

\begin{dfn}
\label{dfn:bridge} Let $N$ be a manifold, $L \subset N$ a link, $F
\subset N$ a \hhs, and denote the complementary handlebodies by $H_1, \
H_2$. The {\it bridge number of $L$ with respect to} $F$ is the minimal number
of arcs in $L' \cap H_1$ for any link $L'$ isotopic to $L$, subject to the
constraint that $L \cap H_1$ and $L \cap H_2$ consists of boundary parallel
arcs.
\end{dfn}

\begin{proof}[Proof of Theorem~\ref{thm:bridge}]

Note that we need to show that the $b(k)$ is bounded above by $g(S)+1$, where
$S$ is an invariant \hhs\ for $M$ found in Theorem \ref{thm:first}.  For a
double cover $f:M \to N$ the {\it singular set } is the set of fixed points.
Similar to Claim~\ref{clm:quotient} we have:

\begin{clm}
\label{clm:singular-arcs}
Let $H$ be a handlebody of genus $g(H)$ and $f$ be an orientation preserving
involution on $H$. Then the singular set of $h$ consists of at most $g(H)+1$
arcs, and these arcs are boundary parallel.
\end{clm}

\begin{rmk}
It is easy to construct involutions that realize the bound above.
\end{rmk}

We prove the claim by induction on $g(H)$.  For balls, this is a result of
Waldhausen~\cite{MR0236916}.  Assume $g(H)>0$.  By the Equivariant Disk
Theorem $H$ admits equivariant disks, either an invariant disk $D$ (Cases One
and Two below) or two disjoint disks $D, \ D'$ that are involutes of each
other (Case Three).  By the  classification of involutions on a disk we know
that in the first case  the intersection of the singular set with $D$ is
either a properly embedded arc (Case One) or in a single point (Case Two).  We
prove the claim in each case:

\medskip

\noindent {\bf Case One:} a single invariant disk $D$ that intersects the
singular set in an arc.  Cutting $H$ along 
$D$ we obtain a handlebody $H_D$.  If $D$ does not separate $H$ we are left with a
genus $g(H)$ handlebody and are done by induction.  If $D$ does separate
$H$, the two complementary pieces are exchanged by $h$ (note that $h|D$ is a
reflection and so orientation reversing) and the singular set of $h$ consists
of a single arc.

\medskip

\noindent {\bf Case Two:} a single invariant disk $D$ that intersects the
singular set in a point.  If $D$ does not separate $H$, cutting $H$ open along
$D$ we obtain a handlebody of genus $g(H) - 1$.  By induction the singular set
consists of at most $g(H)$ boundary parallel arcs. If $D$ separates, cutting
$H$ open along $D$ we obtain two handlebodies (say $H_1$ and $H_2$ of genera
$g_1$ and $g_2$) with $g_1, \ g_2 < g(H)$ and $g_1 + g_2 = g(H)$.  Since
$f|_D$ is orientation preserving, $H_1$ and $H_2$ are invariant under $f$.  By
induction, the singular set in $f|{H_i}$ consists of at most $g_i + 1$
boundary parallel arcs ($1=1,2$). Since $g_1+g_2=g(H)$, adding these numbers
we get $g_1+1 + g_2+1=g(H)+2$.  Luckily, gluing along $D$, two arcs are
identified, becoming a single boundary parallel arc and reducing the number of
singular arcs by one.

\medskip

\noindent {\bf Case Three:} two disjoint disks $D_1$, $D_2$ are exchanged by $f$.
Then $f|_{D_2}, \ f|_{D_2}$ do not admit a fixed point and therefore the
singular set does not intersect $D_1$ or $D_2$.  Cutting $H$ along $D_1$ and
$D_2$, we get at most three components, all handlebodies. If there are one or
two components, the sum of their genera is strictly less than $g(H)$ and by
induction the branch set consists of at most $g(H)$ boundary parallel arcs
that remain boundary parallel after gluing. If there are three components,
then two components (say $H_1$ and $H_2$) are exchanged by the involution and
the last component (say $H_{1,2}$) is invariant. Since $D_1$ and $D_2$ are not
boundary parallel $H_1$ and $H_2$ have positive genus, and therefore the genus
of $H_{1,2}$ is strictly less than $g(H)$. Since $f$ has no fixed points in
$H_1$ or $H_2$ the singular set of $h$ is the same as the singular set of
$h|_{H_{1,2}}$ and the result follows from the inductive hypothesis.  This
completes the proof the Claim~\ref{clm:singular-arcs}.

\medskip

\noindent  Checking the same three cases, one easily proves the
following claim.  To avoid repetition the details are omitted:

\begin{clm}
\label{clm:branch-set-in=handlebody}

Let $H$ be a handlebody and $f:H \to H$ an orientation preserving involution.
Then $H/(f)$ is a handlebody and the branch set of $H$ consists entirely of
boundary parallel arcs.
\end{clm}

Since an involution is injective on the singular set, the branch set of $f:H
\to H/(f)$ has the same number of arcs as its singular set.
Theorem~\ref{thm:bridge} follows from Claims~\ref{clm:singular-arcs} and
\ref{clm:branch-set-in=handlebody}.
\end{proof}

\label{theEnd}

\end{document}